\documentclass[onefignum,onetabnum]{siamart171218}
\usepackage{amsmath, amssymb, bbm, chngcntr}

\usepackage{xr}
\counterwithin{figure}{section}
\usepackage{algorithmic}
\newcommand\algorithmicprocedure{\textbf{procedure}}
\newcommand{\algorithmicendprocedure}{\algorithmicend\ \algorithmicprocedure}
\makeatletter
\newcommand\PROCEDURE[3][default]{%
  \ALC@it
  \algorithmicprocedure\ \textsc{#2}(#3)%
  \ALC@com{#1}%
  \begin{ALC@prc}%
}
\newcommand\ENDPROCEDURE{%
  \end{ALC@prc}%
  \ifthenelse{\boolean{ALC@noend}}{}{%
    \ALC@it\algorithmicendprocedure
  }%
}

\usepackage{booktabs, siunitx}

\newenvironment{ALC@prc}{\begin{ALC@g}}{\end{ALC@g}}
\makeatother

\usepackage{lipsum}
\usepackage{amsfonts}
\usepackage{graphicx}
\usepackage{epstopdf}
\usepackage{algorithmic}
\ifpdf
  \DeclareGraphicsExtensions{.eps,.pdf,.png,.jpg}
\else
  \DeclareGraphicsExtensions{.eps}
\fi

\newsiamremark{remark}{Remark}
\newsiamremark{hypothesis}{Hypothesis}
\crefname{hypothesis}{Hypothesis}{Hypotheses}
\newsiamthm{claim}{Claim}

\headers{A level set approach for circadian phase estimation}{Dae Wook Kim, Minki P. Lee, and Daniel B. Forger}

\title{A Level set Kalman Filter approach to estimate the circadian phase and its uncertainty from wearable data\thanks{Submitted to the editors DATE.
\funding{This work was funded by Human Frontiers Science Program Organization Grant RGP0019/2018 and NSF DMS grant 2052499. DBF is the CSO of Arcascope, a company that makes circadian rhythms software. Both he and the University of Michigan own equity in Arcascope.}}}

\author{Dae Wook Kim\thanks{Department of Mathematics, University of Michigan, Ann Arbor, MI 
  (\email{daewook@umich.edu}, \email{minkilee@umich.edu}).}
\and Minki P. Lee\footnotemark[2]
\and Daniel B. Forger\thanks{Department of Mathematics, University of Michigan, Ann Arbor, MI, and Department of Computational Medicine and Bioinformatics, University of Michigan, Ann Arbor, MI 
  (\email{forger@umich.edu}).}}

\usepackage{amsopn}

\makeatletter
\newcommand*{\addFileDependency}[1]{
  \typeout{(#1)}
  \@addtofilelist{#1}
  \IfFileExists{#1}{}{\typeout{No file #1.}}
}
\makeatother

\newcommand*{\myexternaldocument}[1]{%
    \externaldocument{#1}%
    \addFileDependency{#1.tex}%
    \addFileDependency{#1.aux}%
}

\ifpdf
\hypersetup{
  pdftitle={A Level set Kalman Filter approach to estimate the circadian phase and its uncertainty from wearable data},
  pdfauthor={Dae Wook Kim, and Minki P. Lee, Daniel B. Forger}
}
\fi

\myexternaldocument{ex_supplement}

\begin{document}

\maketitle

\begin{abstract} The circadian clock is an internal timer that coordinates the daily rhythms of behavior and physiology, including sleep and hormone secretion. Accurately tracking the state of the circadian clock, or circadian phase, holds immense potential for precision medicine. Wearable devices present an opportunity to estimate the circadian phase in the real world, as they can non-invasively monitor various physiological outputs influenced by the circadian clock. However, accurately estimating circadian phase from wearable data remains challenging, primarily due to the lack of methods that integrate minute-by-minute wearable data with prior knowledge of the circadian phase. To address this issue, we propose a framework that integrates multi-time scale physiological data to estimate the circadian phase, along with an efficient implementation algorithm based on Bayesian inference and a new state space estimation method called the level set Kalman filter. Our numerical experiments indicate that our approach outperforms previous methods for circadian phase estimation consistently. Furthermore, our method enables us to examine the contribution of noise from different sources to the estimation, which was not feasible with prior methods. We found that internal noise unrelated to external stimuli is a crucial factor in determining estimation results. Lastly, we developed a user-friendly computational package and applied it to real-world data to demonstrate the potential value of our approach. Our results provide a foundation for systematically understanding the real-world dynamics of the circadian clock.
\end{abstract}



\section{Introduction}
Circadian rhythms are $\sim$24hr oscillations in behavioral and physiological processes observed in nearly all living organisms \cite{dibner2010mammalian}. In mammals, including humans, the daily rhythms are coordinated with the external light-dark cycles by an innate timing system, the circadian clock \cite{dibner2010mammalian}. The failure of synchrony between cellular clocks can occur, for instance, due to an alteration of the external environment. Notably, around $80\%$ of the population appears to live a shift work lifestyle \cite{sulli2018training}. This increases the risk for various chronic diseases, such as sleep disorders, psychiatric disorders, cancer, and diabetes \cite{zhu2012circadian}. Thus, tools to calculate sleep schedules rapidly restoring desynchronization have recently received attention \cite{hong2021personalized, walch2016global}. Targets of more than $80\%$ of currently approved drugs have daily rhythmic activity \cite{ruben2019dosing, zhang2014circadian}. As a result, the efficacy and toxicity of diverse drugs, including around 50 anticancer drugs, largely change upon dosing time \cite{kim2020wearable}. Thus, many clinical trials have been performed to develop a pharmaceutical intervention that considers the patient's circadian phase, so-called chronotherapy \cite{panda2019arrival, ruben2019dosing}. To realize chronotherapy, we need the precise measurement of an individual's circadian phase in real-world settings.\\
\indent One of the promising tools for this is wearable technology. Wearable devices already owned by millions of individuals track physiological proxies of clocks, such as rest-activity rhythms and heart rate (HR) outside the laboratory \cite{kim2020wearable}. However, because only the downstream signals can be monitored by wearables, suitable follow-up analysis is required to infer the unobservable state space of the molecular clocks in tissues. One approach is to estimate the phase using differential equation models of the human circadian clock, coupled with activity or light levels recorded from wearables \cite{huang2021predicting, stone2019application, stone2020computational, woelders2017daily}. Specifically, a van der Pol limit cycle model that takes light or activity measurements as the direct input has been proposed to estimate the circadian phase \cite{forger1999simpler, huang2021predicting}. Another recent approach is to demask noisy measurements of the HR rhythm affected by various confounding factors, such as activity, stress, and hormones, using a Bayesian framework with harmonic-regression-plus-autoregressive-noise models and extract a peripheral circadian rhythm in HR \cite{bowman2021method, mayer2022consumer}. This shows the usefulness of the combined approach of wearable technology and mathematical analysis. However, accurate estimation of the circadian phase from wearable data remains far from complete because the information available from wearable monitoring is not fully exploited. For instance, minute-by-minute measurements of rhythmic outputs (e.g., HR) are not currently utilized when estimating the molecular clock phase due to the absence of suitable data assimilation techniques \cite{huang2021predicting, stone2019application, woelders2017daily, zuurbier2015fragmentation}. Moreover, the contribution of noise from different systems to the phase uncertainty remains to be elucidated \cite{huang2021predicting, stone2019application, woelders2017daily}.\\
\indent To circumvent this, a filtering approach that assimilates the knowledge of the system and wearable measurements seems promising \cite{sarkka2013bayesian}. However, conventional data assimilation frameworks are not directly applicable due to the time-scale difference between system estimates and wearable measurements. Specifically, our goal is to estimate the circadian phase, which is a specific time of \textit{day}. However, available information to update the phase estimate is \textit{minute-by-minute} wearable measurements and the clock state is continuously propagated forward \textit{every time}. Thus, straightforward filtering methods for problems where the measurements, the internal system dynamics, and the estimate of interest are on a similar time scale cannot be used to estimate the circadian phase from wearable data.\\
\indent Here, we proposed a generalizable approach for estimating the state of intracellular systems using wearable data. We have applied this approach to assimilate multi-time scale physiological information and estimate the circadian phase, along with its uncertainty, from wearable data. It first extracts daily physiological parameters from short time-scale wearable data using a Bayesian approach. In the second step, the extracted long time-scale information is efficiently and accurately assimilated to estimate the state space of the molecular clocks in tissues using a new state space estimation method called the level set Kalman filter (LSKF). This can account for the contribution of noise from different systems to the estimation, which is impossible with previous methods \cite{huang2021predicting, stone2019application, stone2020computational, woelders2017daily}. Numerical experiments show that our method has a consistent performance improvement over the previous methods. We also apply the method to real-world data to further demonstrate its usefulness.\\
\indent The rest of the paper is structured as follows: In section \ref{sec::2}, we review previous approaches to estimate the circadian phase and describe the novelty of our approach. In section \ref{sec::3}, we define our filtering problem for estimating the circadian phase and describe our approach to solve the problem. Then, our method is tested on various $in$-$silico$ data while changing the noise parameter values, and its performance is compared with the previous approaches in section \ref{sec::numerical simulations}. In section \ref{sec::real data}, our method is applied to real-world wearable data, which shows its usefulness in real-world settings. In section \ref{sec::conclusion}, we conclude the paper by describing the broad applicability of the method, its limitations, and potential future work.
\section{Background: Previous approaches and the novelty of our approach}
\label{sec::2}
In this section, we summarize how our approach to estimate the clock phase is different from previous approaches. Early work done by Brown and colleagues estimated the phase of body temperature circadian rhythms by fitting dynamic (e.g., limit cycle) or simple harmonic models to body temperature data collected in carefully controlled laboratory settings using Bayesian methods \cite{brown2000statistical, brown1992statistical, brown1999statistical, brown19949}. These techniques were used to analyze much of the initial data on human peripheral circadian rhythms. When estimating the central clock phase, many studies perform experiments, such as measuring the time that melatonin collected in dim light crosses a threshold \cite{kim2020wearable, wichniak2017treatment}. However, such methods can only be available when laboratory data are given. Moreover, the methods for the central clock phase estimation are not statistical in that they do not calculate confidence intervals for each measurement.\\
\indent To circumvent this, frameworks to extract the circadian phase from real-world wearable data have been recently proposed \cite{bowman2021method, huang2021predicting, stone2020computational}. The most validated method is to use the limit-cycle oscillator models of the human circadian pacemaker that take in inputs of wearable light measurements or their estimates from wearable activity measurements \cite{forger1999simpler, huang2021predicting, jewett1999revised, kronauer1990quantitative, stone2020computational}. By simulating the models, the estimate of the circadian phase (e.g., the time of the circadian signal minimum) can be computed. Despite this progress, circadian phase estimation in real-life settings remains a challenge. The existing methods do not assimilate the given wearable data despite the potential of data assimilation to improve performance. This is due to the absence of frameworks to combine different sources of physiological information with different time scales (e.g., the molecular clock state and the minute-by-minute HR measurements).\\
\indent To address this, we first proposed a two-step filtering approach for the assimilation of multi-time scale physiological information. It extracts the long time-scale information (i.e., circadian parameters) from short time-scale wearable data. Then, the extracted information is exploited to update the estimate of the unobservable molecular clock state that is continuously propagated. We then developed a numerical algorithm for accurate and efficient implementation of the two-step framework by integrating the LSKF with the Bayesian inference method. The algorithm is easy-to-use and user-friendly so that a broad range of scientists, including physicians, can exploit it for their research. We also provide publicly available user-friendly computer codes to easily implement our method. Lastly, by solving the two-step filtering problem with the proposed algorithm, we analyzed the contribution of the biological noise from different systems to the phase estimation. This analysis is impossible with previous methods \cite{huang2021predicting,stone2019application,stone2020computational, woelders2017daily}, which demonstrates another benefit of our method in addition to performance improvement.
\section{Data assimilation framework to estimate the circadian phase and its uncertainty}
\label{sec::3}
Here, we propose a framework to estimate the circadian phase and its uncertainty based on the LSKF \cite{wang2021level} and the Bayesian inference method \cite{bowman2021method}. This section is organized as follows: In section \ref{sec::3.1}, we explain the mathematical models used to define the process equation and the measurement equation of our filtering framework. In section \ref{sec::3.2}, the process equation and the measurement equation of our LSKF framework are explained. In section \ref{sec::3.3} and section \ref{sec::3.4}, its time-update and measurement-update steps are described, respectively. In section \ref{sec::3.5}, we describe how the circadian clock state estimated on the variable domain can be transformed into that on the time domain. In section \ref{sec::3.6}, we summarize our method and show an example of its application to clearly explain how it works.

\subsection{Mathematical models of the human circadian pacemaker}
\label{sec::3.1}
\subsubsection{Limit-cycle oscillator model of the human circadian pacemaker}
\label{sec::3.1.1}
In this study, we describe our framework with one of the most validated limit-cycle oscillator models of the human circadian clock constructed by Forger, Jewett, and Kronauer \cite{forger1999simpler, huang2021predicting}. The Forger-Jewett-Kronauer (FJK) model can be thought of as a representation of the molecular dynamics of the circadian system. For example, biochemical models of the molecular timekeeping of individual cells can be reduced to the FJK model by averaging on approximate manifolds \cite{forger2002reconciling}. Moreover, a recent Ansatz yields the model from large networks of coupled oscillators \cite{hannay2019macroscopic, hannay2018macroscopic}. This also explains process noise which has a molecular origin \cite{forger2003detailed} and noise in the inputs to the model \cite{forger2004starting} as well as the measurement noise considered by Brown and colleagues \cite{brown1992statistical}. Note that other circadian oscillator models taking light or activity measurements as an input can also be used to implement our filtering framework described below.\\
\indent The human circadian pacemaker affected by light can be described by a van der Pol type oscillator model \cite{forger1999simpler}:

\begin{equation}
    \frac{dx}{dt} = \frac{\pi}{12} (x_c+B(x,x_c,n))
    \label{eq::dxdt}
\end{equation}
\begin{equation}
    \frac{dx_c}{dt} = \frac{\pi}{12} \left[ \mu \Big( x_c - \frac{4}{3}x_c^3 \Big) - x \left\{ \Big( \frac{24}{0.99669\tau_x} \Big)^2 + kB(x,x_c,n) \right\} \right]
    \label{eq::dxcdt}
\end{equation}
where $\mu = 0.23, \tau_x = 24.2$, and $k = 0.55$. The solution trajectory of \cref{eq::dxdt} and \cref{eq::dxcdt} is a limit-cycle oscillation representing the endogenous rhythm of the circadian pacemaker. This oscillator is denoted by \textit{Process P}. $B(x,x_c,n)$ in \cref{eq::dxdt} and \cref{eq::dxcdt}, representing the effect of the photic drive on the circadian pacemaker, and the conversion of light to the photic drive, denoted by \textit{Process L}, are modeled as described below. Light activates photoreceptor activator elements in a ``ready'' state (fraction $1-n$) and converts them to a ``used'' state (fraction $n$) at a rate of $\alpha$, which depends on light intensity $I$ \cref{eq::light}.
\begin{equation}
    \alpha(I) = \alpha_0 \Big( \frac{I}{I_0} \Big)^p
    \label{eq::light}
\end{equation}
where $\alpha_0 = 0.16, p = 0.6$, and $I_0 = 9500$. The used elements are recycled back into the ready state at a rate of $\beta = 0.013$. This is given by 
\begin{equation}
    \frac{dn}{dt} = 60 \Big( \alpha(I)(1-n) - \beta n \Big).
    \label{eq::dndt}
\end{equation}
\begin{figure}[hbt!]
    \label{fig::3.1}
    \centering
    \includegraphics[scale = 0.72]{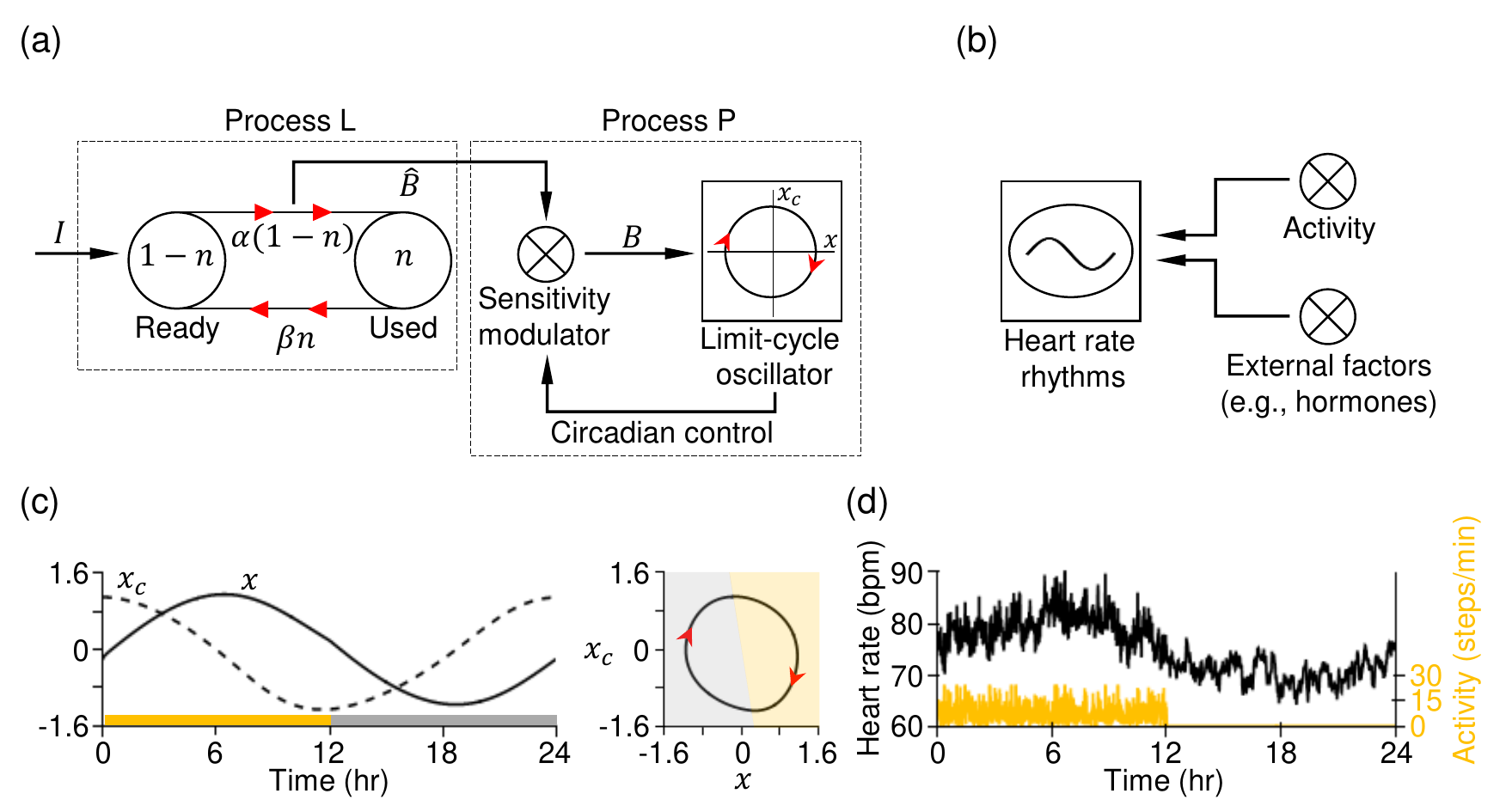}
    \caption{Mathematical models simulating the human circadian timekeeping system. (a,b) Schematic diagram describing the model of the human molecular clocks \cref{eq::dxdt}-\cref{eq::B} (a) and the model of the circadian rhythm in the heart (b) \cref{eq::HR}-\cref{eq::v_t}. (c,d) The exemplary trajectories of $x$ and $x_c$ (c), and $HR_t$ (d). In (c), the light intensity $I(t)=10000$ lux for $t\in\left(0,12 \right)$ and $I(t)=0$ lux for $t\in\left(12,24 \right)$. In (d), the parameter values were adopted from \cite{bowman2021method}.}
    \label{fig:figure3.1.}
\end{figure}
As the elements are activated, it generates a drive onto the circadian pacemaker $\hat{B}$, which is proportional to element flux rate $\alpha(I)(1-n)$ \cref{eq::Bhat}.
\begin{equation}
    \hat{B}(n) = G\alpha(I)(1-n).
    \label{eq::Bhat}
\end{equation}
where $G = 19.875$. The sensitivity of the circadian pacemaker to the drive is modulated in a circadian manner with the state variables $x$ and $x_c$ \cref{eq::B}.   
\begin{equation}
    B(x,x_c,n) = \hat{B}(n)(1-0.4x)(1-0.4x^c).
    \label{eq::B}
\end{equation}
The model diagram of the circadian clock and example solution trajectories of $x$ and $x_c$ are shown in Figures \ref{fig::3.1}a and \ref{fig::3.1}c, respectively.

\subsubsection{Harmonic-regression-plus-correlated-noise model of the HR circadian rhythm}
\label{sec::3.1.2}
The circadian rhythm in HR can be described by a harmonic-regression-plus-first-order-autoregressive model with a linear coefficient that corresponds to how much HR increases per activity (step counts) \cite{bowman2021method, mayer2022consumer}. Specifically, HR oscillates with an approximately 24 hours period \cite{massin2000circadian}, and it increases from this baseline proportionate to activity, matching existing data \cite{scheer2010impact}. Moreover, many external factors, such as cortisol and other hormones \cite{becker2019time}, stress \cite{buckert2014acute}, and caffeine intake \cite{lovallo2006cortisol}, affect HR on the hour timescale. This yields a final model for HR at time $t$, measured in hours:
\begin{equation}
    HR_t = \mu - a \cos\Big( \frac{2\pi}{\tau}(t - \phi^{HR}) \Big) + d \cdot Activity_t + v_t
    \label{eq::HR}
\end{equation}
where
\begin{equation}
    v_t = \alpha \cdot v_{t - \Delta t} + \epsilon_t
    \label{eq::v_t},
\end{equation}
$\mu$ is the basal HR in beats per minute (BPM), $a$ denotes the circadian amplitude, $\phi^{HR}$, representing the time of the circadian HR minimum, denotes the circadian phase as done in \cite{bowman2021method}, the circadian period $\tau = 24$ in hours, $d$ is the increase in HR per unit activity, $|\alpha| < 1$, and the $\epsilon_t$ values are Gaussian random variables with mean zero and variance $\sigma_\epsilon^2$. Note that the noise $v_t$ follows the first-order autoregressive noise process \cref{eq::v_t}: the noise at time $t$ carries over a fraction $\alpha$ of the noise at time $t-\Delta t$. This describes the ongoing effects of extrinsic factors on HR. The independent Gaussian noise $\epsilon_t$ represents new extrinsic effects and measurement error. The model diagram of the HR rhythm and an example trajectory of $HR_t$ are shown in Figures \ref{fig::3.1}b and \ref{fig::3.1}d, respectively.

\subsection{Problem formulation}
\label{sec::3.2}
We consider a continuous-discrete nonlinear filtering model of the circadian timekeeping system whose process equation is formulated with \cref{eq::dxdt}-\cref{eq::B}:

\begin{equation}
    d \mathbf{x}_t = \mathbf{v}(\mathbf{x}_t) dt + \sqrt{\mathbf{K}} dW_t
    \label{eq::dXt}
\end{equation}
where 
\begin{footnotesize}
\[ \mathbf{x}_t = \begin{pmatrix} 
x(t) \\ x_c(t) \\ n(t)
\end{pmatrix},  \quad \mathbf{v}(\mathbf{x}_t) = \begin{pmatrix} 
\frac{\pi}{12}(x_c(t) + B(x,x_c,n)) \\
\frac{\pi}{12} \left[ \mu \Big( x_c(t) - \frac{4}{3}x_c^3(t) \Big) - x(t) \left\{ \Big( \frac{24}{0.99669\tau_x} \Big)^2 + kB(x,x_c,n) \right\} \right] \\
 60 \Big( \alpha(I)(1-n(t)) - \beta n(t) \Big)
\end{pmatrix}, \]
\end{footnotesize}
$W_t$ is a standard 3-dimensional Brownian motion, and $\mathbf{K}\in \mathbb{R}^{3 \times 3}$ is a positive semi-definite continuous process noise matrix, with its decomposition $\mathbf{K} = \sqrt{\mathbf{K}} \sqrt{\mathbf{K}}^T$. The state of the molecular clocks is indirectly measured from the phase of the peripheral HR clock (i.e., $\phi^{HR}$ in \cref{eq::HR}):
\begin{equation}
    \phi_i^{HR} = \phi_i + \phi_{ref} + \epsilon^{\phi}_{i}
    \label{eq::HRphase}
\end{equation}
where $\phi_i^{HR}$ denotes the HR phase on day $i$, $\phi_i$ is the phase of the circadian pacemaker on day $i$, $\phi_{ref}=-1$ denotes the average difference between the two phases reported in the previous studies \cite{bowman2021method, forger1999simpler, wichniak2017treatment}, and $\epsilon^{\phi}_{i}$ is the zero-mean Gaussian measurement noise on day $i$. For a given day $i$, we define $\phi_i$ to be 
\begin{equation}
    \phi_i = \{ t_i \pmod{24} : \mathbb{E}[x(t_i)] \leq \mathbb{E}[x(t)], \forall t \in [24(i-1),24i] \}.
    \label{eq::internal phase}
\end{equation}
Note that $\mathbb{E}[x(t)]$ denotes the expectation of $x(t)$ with respect to the probability distribution of $x(t) = \left[x(t), x_{c}(t), n(t)\right]'$ in \cref{eq::dXt}. The mean dynamics computed using the averaged velocity level set time-update method (see section \ref{sec::3.3.1} below) has a limit cycle of period $24$hr, and there is only one minimum point of $x(t)$ for each cycle (Supplementary Materials and \cite{diekman2018reentrainment}).\\
\indent The measurement equation of our filtering framework is represented by the relationship between the two clock states \cref{eq::HRphase}-\cref{eq::internal phase}. The method used to obtain the HR phase estimate from wearable data \cite{bowman2021method} is described in section \ref{sec::3.4}.

\subsection{Time-update step}
\label{sec::3.3}
Here, we use the time-update method of the LSKF to propagate the clock state estimate forward in time because it better fits our problem settings, compared to other flavors of Kalman filter, such as the continuous-discrete cubature Kalman filter (CD-CKF). Specifically, although the time-update of the CD-CKF (with Ito-Taylor expansion of order 1.5) and that of the LSKF converge to the same limit, the LSKF method can achieve a higher order of convergence \cite{wang2021level}. Moreover, the LSKF can outperform the CD-CKF with sufficient time-step subdivisions, which is the choice for challenging filtering problems \cite{arasaratnam2010cubature, kulikov2017accurate}, when the measurement interval is large. Thus, the LSKF is suitable for our data assimilation problem with a significantly large measurement interval: There is only one measurement (i.e., $\phi^{HR}_{i}$) on each day. Another possible time-update method is to use particle filters \cite{mott2011model}. However, their computational costs are prohibitive in many real-world applications because they typically require the simulation of individual stochastic trajectories. Unlike this, the LSKF computes the propagation of the estimate in time efficiently by solving coupled ordinary differential equations for the Gaussian level set.\\ 
\indent Another reason to select the LSKF time-update method is that it is user-friendly and easy-to-use. For example, while the CD-CKF requires the spatial partial derivatives of the drift function and user-defined time-step subdivisions explicitly, the LSKF only requires knowledge of the drift function \cite{wang2021level}. Moreover, the LSKF time-update step can be easily performed using any adaptive ODE solver without an appropriate time discretization.\\
\indent In section \ref{sec::3.3.1}, we review the time-update method of the LSKF. In section  \ref{sec::3.3.2}, we describe how the best possible estimate of $\mathbf{x}$ in \cref{eq::dXt} can be propagated using the time-update step of the LSKF.

\subsubsection{Review of the time-update method of the LSKF}
\label{sec::3.3.1}
The probability density $u(\mathbf{x},t)$ of $\mathbf{x}$ in \cref{eq::dXt} satisfies the Fokker-Plank equation in Ito's calculus \cite{bressloff2014stochastic}:
\begin{equation}
    \frac{\partial u(\mathbf{x},t)}{\partial t} = \frac{1}{2} \nabla \cdot \mathbf{K} \nabla u(\mathbf{x},t) - \nabla \cdot (\mathbf{v}(\mathbf{x})u(\mathbf{x},t))
    \label{eq::fokkerplank}.
\end{equation}
Without loss of generality (WLOG), we may assume the drift function at $\mathbf{x} = 0$ is zero (i.e., $\mathbf{v}(0) = 0$). Then we approximate \cref{eq::fokkerplank} by taking a linear approximation of $\mathbf{v}$ (i.e., $\mathbf{v}(\mathbf{x}) \approx \mathbf{Jx}$, where $\mathbf{J}$ is the Jacobian matrix):
\begin{equation}
    \frac{\partial u(\mathbf{x},t)}{\partial t} = \frac{1}{2} \nabla \cdot \mathbf{K} \nabla u(\mathbf{x},t) - \nabla \cdot (\mathbf{Jx}u(\mathbf{x},t)).
    \label{eq::revised_fokkerplank}
\end{equation}
The auxiliary function $F(\mathbf{x},t):=u(\mathbf{x},t)/u(\mathbf{0},t)$ is introduced to define a level set of a Gaussian distribution as 
\begin{equation}
    \mathcal{L}(F(\mathbf{x},t),c) := \{ \mathbf{x} : F(\mathbf{x},t) = c \}
    \label{eq::levelset}
\end{equation}
where $c \in (0,1)$ is some fixed scalar constant. Note that $\mathcal{L}(F(\mathbf{x},t),c)$ is an ellipsoid if $u(\mathbf{x},0)$ is given by a Gaussian function
\begin{equation}
    u(\mathbf{x},0) = \frac{1}{\sqrt{(2\pi)^d \det \Sigma}} \exp \Big( \frac{-\mathbf{x}^T \Sigma^{-1} \mathbf{x}}{2} \Big)
    \label{eq::levelset_ellipsoid}
\end{equation}
where $d$ is a dimension of $\mathbf{x}$, $\Sigma$ is the covariance matrix, and it is assumed WLOG that $\mathbf{x}$ is centered at 0 at time 0. Since $F$ varies over time, $\mathcal{L}(t)$ propagates in space. To describe its propagation, we consider a velocity of level set $\mathbf{v}_{\mathcal{L}}$ defined by a velocity field satisfying the level-set equation \cite{sethian1985curvature}:
\begin{equation}
    \frac{\partial F}{\partial t} + \mathbf{v}_{\mathcal{L}} \cdot \nabla F = 0.
    \label{eq::sethian}
\end{equation}
A velocity of the level set for $\mathcal{L}(0)$ defined in \cref{eq::levelset} can be explicitly calculated if  $u(\mathbf{x},0)$ is given by a Gaussian function in \cref{eq::levelset_ellipsoid} \cite{wang2021level} as follows:
\begin{equation}
    \mathbf{v}_{\mathcal{L}} = \mathbf{J}\mathbf{x} + \frac{1}{2} \mathbf{K} \Sigma^{-1} \mathbf{x}. 
    \label{eq::velocity}    
\end{equation}
Using \cref{eq::velocity}, we can show that the solution of a Fokker-Plank equation with linear drift function is a Gaussian distribution if an initial condition $u(\mathbf{x},0)$ is given by a Gaussian distribution as follows: 
\begin{corollary} \label{cor1}
  A Fokker-Plank equation preserves a Gaussian distribution with a linear drift function \cref{eq::revised_fokkerplank}
  \begin{proof}
    
    \cref{eq::velocity} is a linear transformation (i.e., a linear velocity field) that instantaneously propagates every level set, independent of the scaling constant $c$. Since a linear transformation maps Gaussian to Gaussian, the initial Gaussian distribution is mapped to a Gaussian under \cref{eq::velocity}. Hence, $u(\mathbf{x},t)$, the propagated distribution at time $t$, is also a Gaussian distribution, and the velocity field $\mathbf{v}_{\mathcal{L}}$ is always well-defined for all $t$.
  \end{proof}
\end{corollary}
Note that the proof of corollary \ref{cor1} is adopted from \cite{wang2021level} and corollary \ref{cor1} can be considered as a corollary of Equation (29) in  \cite{kalman1961new}.\\
\indent We now describe the numerical algorithm for the time-update step. Inspired from corollary \ref{cor1}, we can consider that propagation of a solution of the Fokker-Plank equation with linear drift function for time $t$ \cref{eq::revised_fokkerplank} is equivalent to tracking its ellipsoid level sets. Hence, we track the movement of the column vectors of $\mathbf{M}$ where $\Sigma = \mathbf{M} \mathbf{M}^T$ since they represent the unique level set for the time update \cite{arasaratnam2009cubature}. Specifically, let $\Sigma(0) = \mathbf{M}(0) \mathbf{M}(0)^T$ be the covariance matrix of an initial condition $u(\mathbf{x},0)$, and $\mathbf{M}(0) = \begin{bmatrix} \mathbf{x}_1(0) & \cdots & \mathbf{x}_d(0) \end{bmatrix}$ where $\mathbf{x}_i(0)$ is the $i$-th column vector $\mathbf{x}_i$ of $\mathbf{M}(0)$. Then, we can interpret $\mathbf{x}_i(0)$ as a representative point of one of level sets at time $t = 0$ \cite{wang2021level}. It travels at the speed defined by a velocity of level set in \cref{eq::velocity}. Let the position of the propagated points at time $t$ be denoted by $\mathbf{M}(t) = \begin{bmatrix} \mathbf{x}_1(t) & \cdots & \mathbf{x}_d(t) \end{bmatrix}$. Then, $\Sigma(t) = \mathbf{M}(t) \mathbf{M}(t)^T$ is the covariance matrix for the Gaussian that is a solution of \cref{eq::revised_fokkerplank} at time $t$ because the velocity field defined in \cref{eq::velocity} is linear and a Gaussian is preserved by a linear transformation. By substituting each column vector of $\mathbf{M}(t)$ for $\mathbf{x}$ in \cref{eq::velocity} and approximating the Jacobian $\mathbf{J}$ with the central difference in velocity, we obtain
\begin{equation}
    \frac{d\mathbf{x}_i(t)}{dt} = \mathbf{v}(\bar{\mathbf{x}} + \mathbf{x}_i) - \frac{1}{2d} \sum_{i=1}^d (\mathbf{v}(\bar{\mathbf{x}} + \mathbf{x}_i) + \mathbf{v}(\bar{\mathbf{x}} - \mathbf{x}_i)) + \frac{1}{2} \mathbf{K}(\mathbf{M}(t)^T)^{-1} e_i
    \label{eq::avg-difference}
\end{equation}
where $\bar{\mathbf{x}}$ is the mean of the Gaussian and $e_i$ is the $i$th unit vector with all entries 0 except for $i$th entry. The matrix form of \cref{eq::avg-difference} is 
\begin{equation}
    \frac{d\mathbf{M}}{dt} = \mathbf{v}(\bar{\mathbf{x}} + \mathbf{M}) - \frac{1}{2d} \sum_{i=1}^d (\mathbf{v}(\bar{\mathbf{x}} + \mathbf{x}_i) + \mathbf{v}(\bar{\mathbf{x}} - \mathbf{x}_i)) + \frac{1}{2} \mathbf{K}(\mathbf{M}^T)^{-1}.
    \label{eq::avg-difference2}
\end{equation}
With the averaged velocity (i.e., the second term of \cref{eq::avg-difference}), we can also set the velocity of $\bar{\mathbf{x}}$ as follows:    
\begin{equation}
    \frac{d\bar{\mathbf{x}}}{dt} = \frac{1}{2d} \sum_{i=1}^d (\mathbf{v}(\bar{\mathbf{x}} + \mathbf{x}_i) + \mathbf{v}(\bar{\mathbf{x}} - \mathbf{x}_i)).
    \label{eq::avg-mean-diff}
\end{equation}
We concatenate $\bar{\mathbf{x}}$ and $\mathbf{M}$ as a variable $(\bar{\mathbf{x}}|\mathbf{M})$ of $d \times (d+1)$ dimension and obtain a nonlinear ordinary differential equation (ODE) system. By solving it with any standard ODE solver, we can complete the time-update step between the measurements. 

\subsubsection{Forward propagation of the estimate of the clock state}
\label{sec::3.3.2}
\begin{figure}[hbt!]
    \label{fig::3.2}
    \centering
    \includegraphics[scale = 0.74]{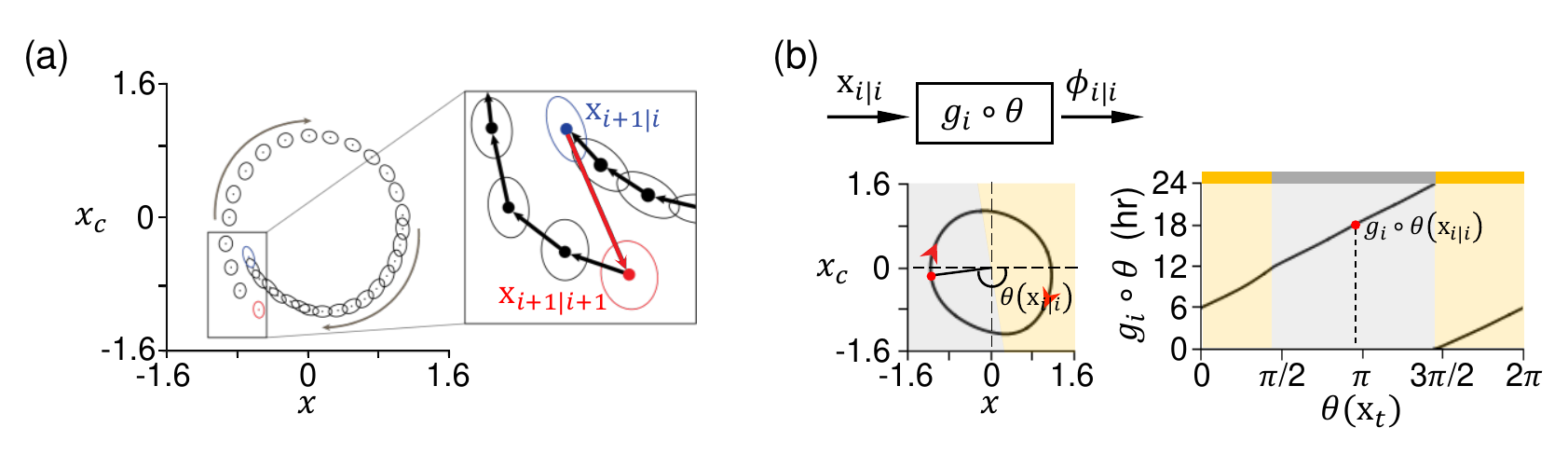}
    \caption{Key steps of our filtering framework. (a) Graphical illustration of the time-update step and the measurement-update step. By the time-update method, the ellipsoid level set of $\mathbf{x}$ represented as a circle is propagated forward (black), and the circadian clock state is predicted (blue). The predicted state is updated by the measurement-update method (red). (b) Transformation of the clock state from the variable domain to the time domain.}
    \label{fig:figure3.2.}
\end{figure}

Here, we apply the LSKF time-update method to our process equation \cref{eq::dXt} to propagate the estimate of $\mathbf{x}$ forward until the subsequent measurement is available. To track the propagation of the level set by solving \cref{eq::avg-difference2} and \cref{eq::avg-mean-diff}, the light intensity $I$ needs to be defined for each time. However, ubiquitous consumer-grade wearable devices typically do not record light levels. To address this problem, Huang et al. showed that light levels are correlated to activity levels in real-life settings and thus light levels can be estimated from activity levels using a carefully constructed steps-to-light function \cite{huang2021predicting}. Thus, the ODE model \cref{eq::dxdt}-\cref{eq::B}, taking in inputs of light levels estimated from activity levels using the function, showed reasonable accuracy in estimating the circadian phase in multiple clinical datasets \cite{huang2021predicting}. Importantly, activity is better at estimating the phase than light in a shift-worker cohort. We adopted this validated steps-to-light piecewise function denoted by $L$ \cref{eq::acttolux} in this study: The step count measured by wearables at time $t$ denoted by $a(t)$ is converted to the light input $I$ with $L$ \cref{eq::acttolux}.
\begin{equation}
L(a(t)) = \begin{cases} 
0 \text{ lux} & \text{if } a(t) \leq 0 \\
100 \text{ lux} & \text{if } a(t) \in (0, 0.1m) \\
200 \text{ lux} & \text{if } a(t) \in [0.1m, 0.25m) \\
500 \text{ lux} & \text{if } a(t) \in [0.25m, 0.4m) \\
2000 \text{ lux} & \text{otherwise}
\end{cases}
\label{eq::acttolux}
\end{equation}
where $m = \frac{\max_{s_1\leq t \leq s2}(a(t))}{2}$, $s_1$ and $s_2$ denote the starting time and ending time of data collection, respectively. By substituting $L(a(t))$ for $I(t)$ in \cref{eq::dXt} and solving \cref{eq::avg-difference2} and \cref{eq::avg-mean-diff}, we can compute the propagation of $\mathbf{x}$ until a new measurement is available, which is shown as black circles in Figure \ref{fig::3.2}a. This allows the prediction of the circadian clock state on day $i+1$ when measurements until day $i$ are given (i.e., $\mathbf{x}_{i+1|i}\sim N(\bar{\mathbf{x}}_{i+1|i},\Sigma_{i+1|i})$), which is shown as a blue circle in Figure \ref{fig::3.2}a.
\subsection{Measurement-update step} \label{sec::3.4}
In section \ref{sec::3.4.1}, we review the measurement\hspace{0.01cm}-update method used in the LSKF. In section \ref{sec::3.4.2}, we describe how the estimate of $\mathbf{x}$ in \cref{eq::dXt} is updated using the method when the subsequent measurement is given.

\subsubsection{Review of the measurement-update of the LSKF} \label{sec::3.4.1}
For the measure-ment-update, the method from the square root continuous-discrete cubature Kalman filter \cite{arasaratnam2010cubature, wang2021level} is used. The method was adopted because it can accommodate a positive semi-definite matrix $\mathbf{M}$ calculated in the LSKF time-update step, which is required for the reliability of the method \cite{arasaratnam2010cubature, wang2021level}. Because our notations are different from those used in \cite{arasaratnam2010cubature, wang2021level}, we restate the measurement-update algorithm in Algorithm \ref{alg::measurement-update}.
\begin{algorithm}[hbt!]
\caption{Measurement-update step}
\label{alg::measurement-update}
\begin{algorithmic}[1]
\REQUIRE Predicted mean after the time-update step $\bar{\mathbf{x}}$, a factorization of the predicted covariance matrix after the time-update step $\mathbf{M}$ such that $\Sigma = \mathbf{M} \mathbf{M}^T$, measurement $\mathbf{z}$, a factorization of measurement noise matrix $\sqrt{\mathbf{R}}$, and a measurement function $\mathbf{h}$
\STATE Find the concatenated matrix of the cubature points
    \begin{equation}
        \mathbf{N} = \mathbf{\Bar{x}} + \sqrt{2d} (\mathbf{M}|-\mathbf{M}),
    \end{equation}
    by applying a vector-matrix addition to $\bar{\mathbf{x}}$ and each column of the concatenated matrix $(\mathbf{M}|-\mathbf{M})$. Note that $\mathbf{N}$ is a $d\times 2d$ matrix if $\bar{\mathbf{x}}$ is a $d\times1$ matrix.\\
\STATE Evaluate the measurement function $\mathbf{h}$ on each column of $\mathbf{N}$ to compute the propagated cubature points as follows:
    \begin{equation}
        \mathbf{Z} = \mathbf{h}(\mathbf{N}).
    \end{equation}
    Note that $\mathbf{Z}$ is $1\times2d$ matrix if $\mathbf{h}$ is a scalar-valued function.
\STATE Estimate the predicted measurement from the propagated cubature points.
    \begin{equation}
        \mathbf{\Bar{z}} = \frac{1}{2d} \sum_{i=1}^{2d} \mathbf{Z}_i,
    \end{equation}
    where $\mathbf{Z}_{i}$ is the $i$th column of $\mathbf{Z}$.
\STATE Perform the QR-factorization to obtain matrix $\mathbf{T_{11}}$ and $\mathbf{T_{21}}$
    \begin{equation}
        \begin{pmatrix} 
        \mathbf{T}_{11} & \mathbf{O} \\
        \mathbf{T}_{21} & \mathbf{T}_{22}
        \end{pmatrix} = \operatorname{qr} \begin{pmatrix} 
        \mathbf{Z} & \sqrt{\mathbf{R}} \\
        \mathbf{N} & \mathbf{O}
        \end{pmatrix},
    \end{equation} 
    where $\mathbf{O}$ is a zero matrix of appropriate size, such that the QR factorization can be performed.
\STATE Use a backward stable solver to compute the cubature gain $\mathbf{W}$ such that
    \begin{equation}
        \mathbf{T}_{21}=\mathbf{W}\mathbf{T}_{11}
    \end{equation}
\STATE Use the computed cubature gain $\mathbf{W}$ to estimate the corrected mean 
    \begin{equation}
        \mathbf{\hat{x}} = \mathbf{\Bar{x}} + \mathbf{W}(\mathbf{z}-\mathbf{\Bar{z}})
    \end{equation}
\STATE Estimate a factorization of the corrected covariance matrix
    \begin{equation}
        \mathbf{\hat{M}} = \mathbf{T}_{22}
    \end{equation}
\RETURN The corrected mean $\mathbf{\hat{x}}$ and a factorization of the corrected covariance matrix $\mathbf{\hat{M}}$
\end{algorithmic}
\end{algorithm}
\subsubsection{Update of the estimate of the circadian clock state with wearable measurements} \label{sec::3.4.2} 
Here, we describe how the predicted estimate of $\mathbf{x}_{i+1|i}$ in section \ref{sec::3.3.2} is updated with wearable data. Because wearable HR data are minute-by-minute measurements, they cannot be directly assimilated into the predicted circadian phase. Thus, we first extract the HR phase on each day $i$ (i.e., $\phi_{i}^{HR}$ in \cref{eq::HRphase}) from wearable HR and activity data. Specifically, we fit \cref{eq::HR} to the HR and activity data using the recently developed Bayesian inference framework \cite{bowman2021method} that is based on Goodman and Weare's affine-invariant MCMC algorithm, which prevents potential bias from large gaps in wearable data \cite{goodman2010ensemble}. This returns the mean and variance of the HR phase estimate computed from its posterior distribution on day $i+1$. We use them as the measurement $\phi_{i+1}^{HR}$ in \cref{eq::HRphase} (i.e., $\mathbf{z}$ in Algorithm \ref{alg::measurement-update}) and the measurement noise $\epsilon_{i+1}^{\phi}$ in \cref{eq::HRphase} (i.e., $\mathbf{R}$ in Algorithm \ref{alg::measurement-update}). This allows implementation of the measurement-update method with the measurement equation \cref{eq::HRphase} denoted by $\mathbf{h}$ in Algorithm \ref{alg::measurement-update}. This updates $\mathbf{x}_{i+1|i}$ predicted by the time-update method and returns the corrected estimate $\mathbf{x}_{i+1|i+1}\sim N(\bar{\mathbf{x}}_{i+1|i+1},\Sigma_{i+1|i+1})$, which is shown as a red circle in Figure \ref{fig::3.2}a. This can be transformed into the estimate $\phi_{i+1|i+1}$ on the time domain (see section \ref{sec::3.5} below for details). 

\subsection{Transformation of the circadian clock state estimate from the variable domain to the time domain} \label{sec::3.5}
Our method estimates the circadian clock state $\mathbf{x}_{i|i}\sim N(\bar{\mathbf{x}}_{i|i},\Sigma_{i|i})$. Here, we describe how the estimated distribution can be transformed from the variable domain to the time domain. We first define a function $\theta(\mathbf{x})$ to establish a relationship between the state vector $\mathbf{x}$ and the angle formed by the projection of $\mathbf{x}$ on the $\left(x,x_c\right)$-plane with the origin as follows:
\begin{equation}
    \label{eq::state_to_angle}
    \theta(\mathbf{x}) = \begin{cases}
    -\arctan(\frac{x_c}{x}) & \text{if } x \geq 0, x_c \leq 0 \\
    \pi - \arctan(\frac{x_c}{x}) & \text{if } x < 0 \\
    2\pi - \arctan(\frac{x_c}{x}) & \text{if } x \geq 0, x_c > 0 
    \end{cases}
\end{equation}
Then, we construct an interpolation function $g_i: [0,2\pi] \rightarrow [t_i, t_{i+1}]$ that relates $\theta(\mathbf{x}_t)$ to time $t$ where $t_i$ and $t_{i+1}$ is the midnight (00:00 am) of day $i$ and day $i+1$, respectively, using the trajectory of $\mathbf{x}$ on day $i$ obtained by propagating it with the level-set method. Finally, the composition $(g_i \circ \theta)(\mathbf{x}_t)$ allows us to relate $\mathbf{x}_t$ to $t$ during day $i$. As a result, we can obtain the circadian clock state on the time domain \cref{eq::transangle} as shown in Figure \ref{fig::3.2}b.  
\begin{equation}
    \label{eq::transangle}
    \phi_{i|i} = (g_i \circ \theta)(\mathbf{x}_{i|i})
\end{equation}
where $\mathbf{x}_{i|i}\sim N(\bar{\mathbf{x}}_{i|i},\Sigma_{i|i})$. To sample from the nonlinearly transformed random variable $\phi_{i|i}$, we used a Monte Carlo approach: we took a large number of $N$ random samples from $\mathbf{x}_{i|i}\sim N(\bar{\mathbf{x}}_{i|i},\Sigma_{i|i})$, and computed their transformed values that are $N$ samples of $\phi_{i|i}$.

\subsection{Summary of the filtering algorithm and an example of its application} \label{sec::3.6}
\begin{figure}[hbt!]
    \label{fig::3.3}
    \centering
    \includegraphics[scale = 0.74]{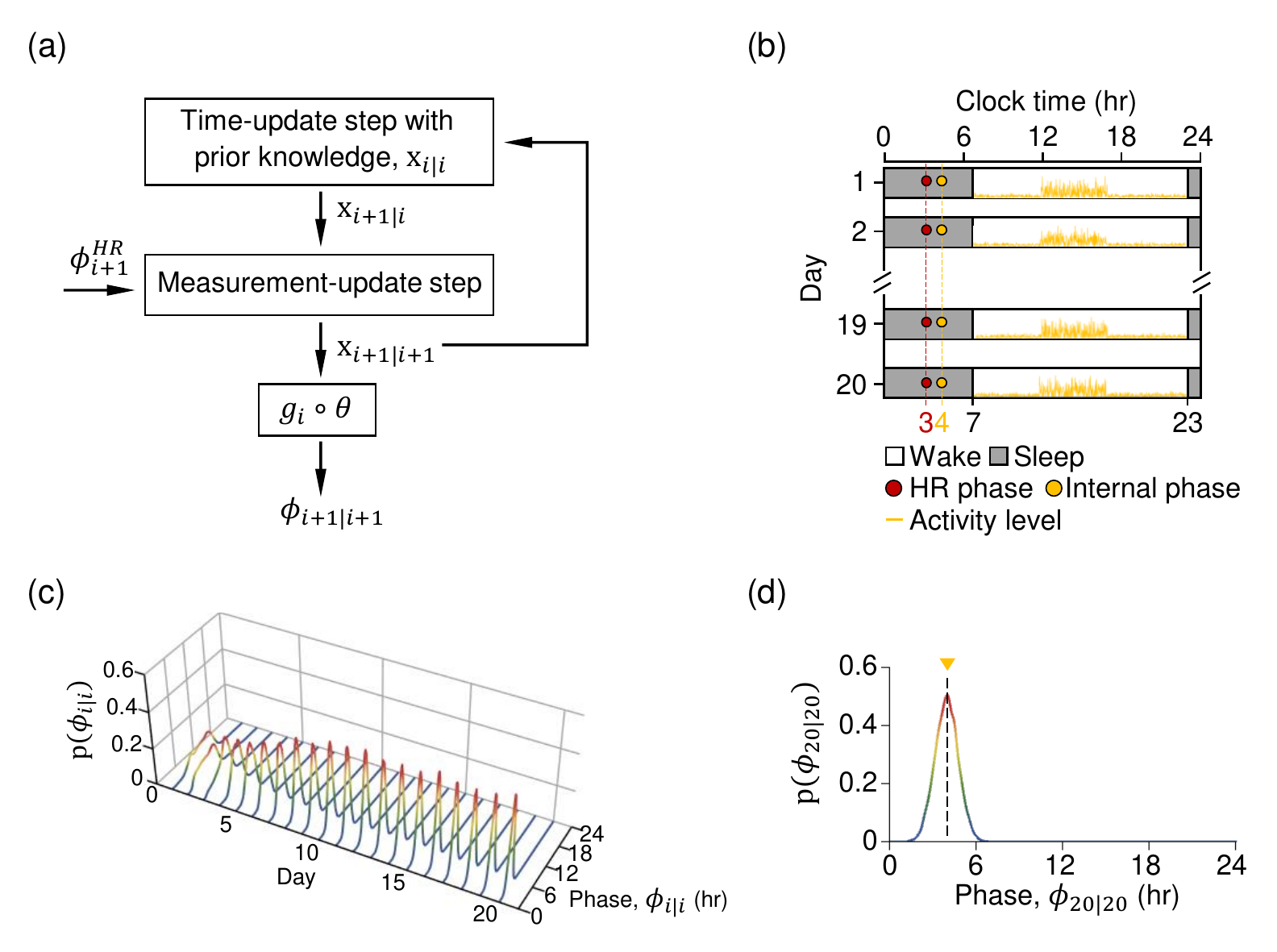}
    \caption{Kalman filtering for the circadian phase estimation. (a) Diagram describing the Kalman filter to estimate the evolution of the posterior distributions of the internal circadian phase $\phi_{i|i}$. (b) A simulated scenario that captures a typical lifestyle. Sleep offset and onset times were set as $7$hr and $23$hr, respectively. The wake and sleep times are represented as white and black box. Activity level $a(t)$ in \cref{eq::acttolux} is denoted by yellow line. The HR phase and the internal phase were set as $3$hr and $4$hr, respectively. Please see Scenario 1 described in section \ref{sec::4.1} and Table \ref{tab::summary_scenario_case} for more details. (c) The evolution of the posterior distributions of $\phi_{i|i}$ over days. Here, a guess of $\bar{\mathbf{x}}_{t_{0}}$ was set as $\left[x(t_0), x_c(t_0), n(t_0) \right]=\left[1, 0, 0.5 \right]$ that is completely different to the true $\left[-0.61, -0.76, 0.34 \right]$. A guess of $\Sigma(t_{0})$ was set as $0.1\cdot I$ and $K=10^{-2}\cdot I$ where $I$ is the identity matrix of order $3$. (d) The posterior distribution on the last day (i.e., $\phi_{20|20}$). The triangle and dashed line represent the true internal phase $\phi_{i}^{true}=4$hr.}
    \label{fig:figure3.3.}
\end{figure}

Our filtering framework can be illustrated in Figure \ref{fig::3.3}a. We predict the estimate on day $i+1$, $\mathbf{x}_{i+1|i}$, from prior knowledge $\mathbf{x}_{i|i}$ as described in section \ref{sec::3.3}. Then, it is updated with $\phi^{HR}_{i+1}$ extracted by applying the Bayesian method \cite{bowman2021method} to wearable data on day $i+1$ as described in section \ref{sec::3.4}. This returns the estimate $\mathbf{x}_{i+1|i+1}$ that is transformed into the estimate on the time domain $\phi_{i+1|i+1}$, as described in section \ref{sec::3.5}. By iterating this procedure, we can estimate the phase on each day in a consecutive manner. This can be summarized as Algorithm \ref{alg::LSKF circadian}.\\
\indent To clearly illustrate how our algorithm works, we applied it to $in$-$silico$ data and showed its estimation results. Specifically, we first created a virtual scenario mimicking a typical human lifestyle: a virtual subject sleeps from 23hr to 7hr and regularly acts during the waking time, as described in Figure \ref{fig::3.3}b. The simulated HR rhythm becomes minimum at the midpoint of sleep (i.e., the true HR phase $=3$hr) following the previous work \cite{bowman2021method}, and the true internal phase of the pacemaker on day $i$ $\phi_{i}^{true}=4$hr. From these activity and HR data, the evolution of posterior distributions of the internal circadian phase $\phi_{i|i}$ can be generated by our method, as shown in Figure \ref{fig::3.3}c. In Figure \ref{fig::3.3}c, a guess of the initial mean state $\bar{\mathbf{x}}_{t_0}$ was completely wrong. Moreover, a guess of the initial covariance matrix of $\mathbf{x}$ $\Sigma(t_{0})$ was set to be large so that initial priors were not highly informative, resulting in the broad distribution on day $1$. Despite this circumstance, our method can accurately estimate the internal circadian phase after a few days, as shown in Figure \ref{fig::3.3}c and \ref{fig::3.3}d. 

\begin{algorithm}[hbt!]
\caption{Data assimilation for the circadian phase estimation}
\label{alg::LSKF circadian}
\begin{algorithmic}[1]
\REQUIRE A guess of initial mean state $\bar{\mathbf{x}}_{t_0}$ at time $t_0$, a factorization of a guess of covariance matrix $\mathbf{M}(t_{0})$, the HR phase estimates $\phi_1^{HR}, \cdots, \phi_n^{HR}$ and its estimate errors $\epsilon^{\phi}_{1}, \cdots, \epsilon^{\phi}_{n}$, on day $1, \cdots, n$ 
\FOR{day $i$ from day 1 to day $n$}
\STATE Time-update: solve \cref{eq::avg-difference2} and \cref{eq::avg-mean-diff} with initial conditions $\bar{\mathbf{x}}_{i-1|i-1}$ and $\mathbf{M}_{i-1|i-1}$ such that $\Sigma_{i-1|i-1}=\mathbf{M}_{i-1|i-1}\mathbf{M}_{i-1|i-1}^T$ from $t = \phi_{i-1|i-1}$ to $t = \phi_{i|i-1}$ using a numerical ODE solver to obtain $\mathbf{x}_{i|i-1}\sim N(\bar{\mathbf{x}}_{i|i-1},\Sigma_{i|i-1})$.
\STATE Measurement-update: perform Algorithm \ref{alg::measurement-update} with $\bar{\mathbf{x}}_{i|i-1}$, $\mathbf{M}_{i|i-1}$, $\phi_i^{HR}$, and $\sqrt{\epsilon^{\phi}_{i}}$, and the measurement function \cref{eq::HRphase} to compute $\mathbf{x}_{i|i}\sim N(\bar{\mathbf{x}}_{i|i},\Sigma_{i|i})$.
\STATE Transformation into the time domain: $\phi_{i|i}=(g_i \circ \theta)(\mathbf{x}_{i|i})$ where $\mathbf{x}_{i|i}\sim N(\bar{\mathbf{x}}_{i|i},\Sigma_{i|i})$ 
\RETURN $\phi_{i|i}$
\ENDFOR
\end{algorithmic}
\end{algorithm}
\section{Numerical study} 
\label{sec::numerical simulations}
Here, we perform a series of numerical experiments to study the usefulness of our method. This section is structured as follows: In section \ref{sec::4.1}, we investigate the relationship between the process noise $\mathbf{K}$ in the circadian clock and the phase estimate, which is impossible with previous methods \cite{huang2021predicting, stone2019application, woelders2017daily}. In section \ref{sec::4.2}, we show that our method has an overall performance improvement over previous methods \cite{huang2021predicting,bowman2021method}.

\subsection{Relationship between the process noise and the phase estimate} 
\label{sec::4.1}
\begin{figure}[hbt!]
    \label{fig::4.1}
    \centering
    \includegraphics[scale = 0.665]{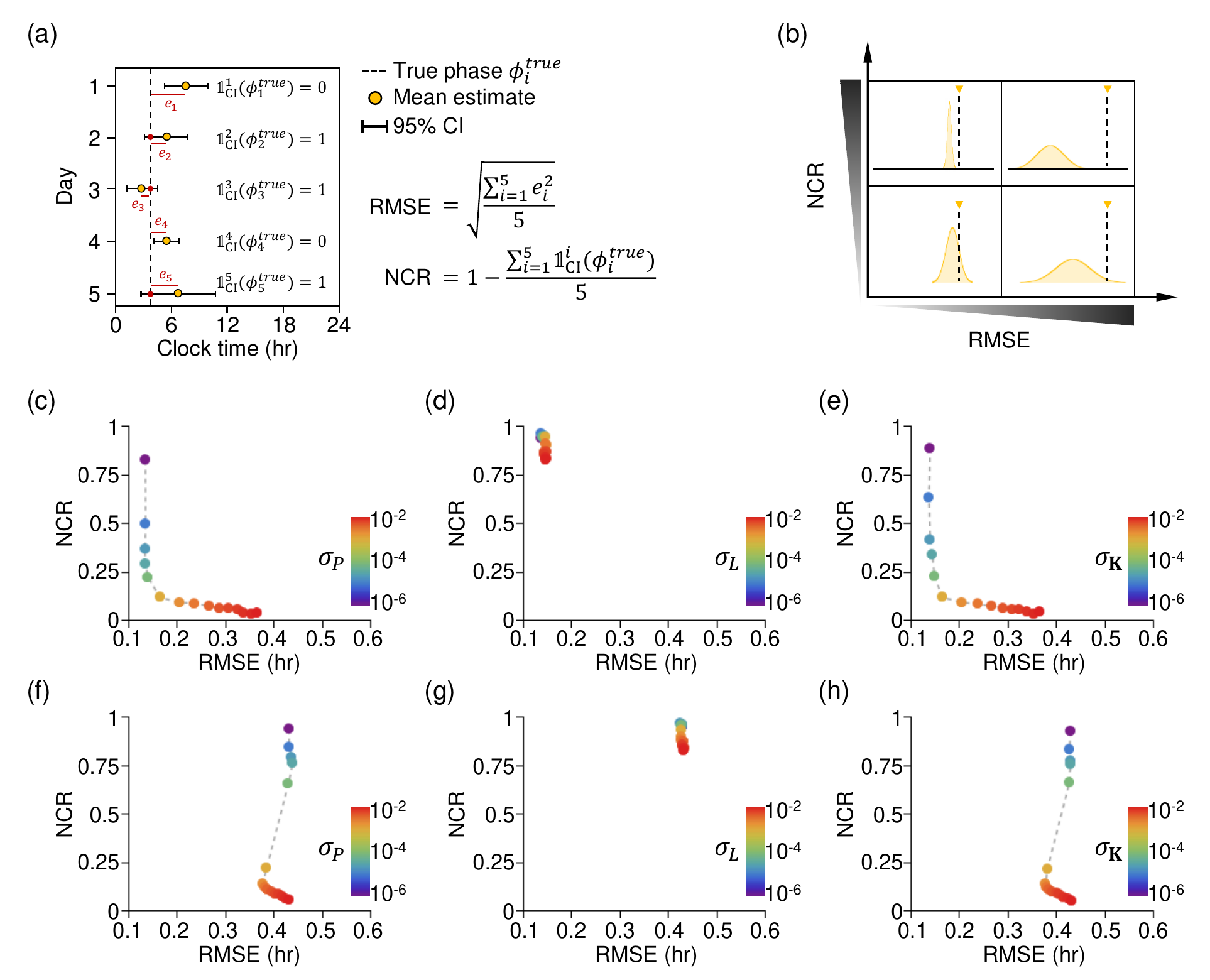}
    \caption{The association of the circadian phase estimate with the noise in the circadian system. (a) The schematic illustration of the definition of the RMSE and the NCR. (b) The relationship of the phase estimate with the RMSE and the NCR. (c-h) The change of the RMSE and the NCR upon the noise magnitude in Scenario 1 (c-e) and in Scenario 2 (f-h). In (c) and (f), only the noise in \textit{Process P} (i.e., $\sigma_{P}$) was considered. In (d) and (g), only the noise in \textit{Process L} (i.e., $\sigma_{L}$) was considered. In (e) and (h), the two noises were considered together.}
    \label{fig:figure4.1.}
\end{figure}

\begin{figure}[hbt!]
    \label{fig::4.2}
    \centering
    \includegraphics[scale = 0.665]{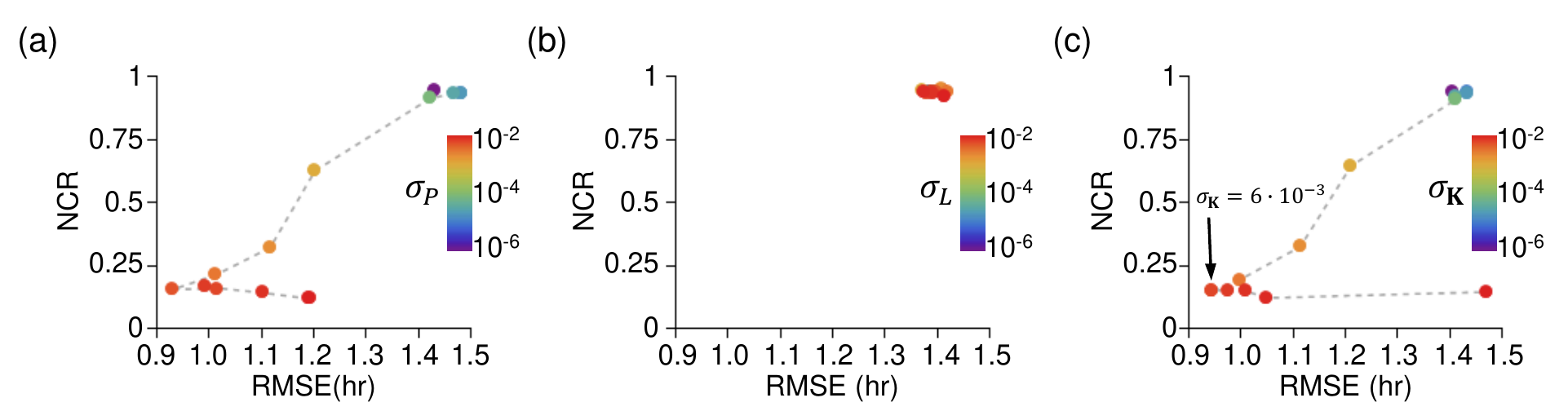}
    \caption{The association of the circadian phase estimate with the noise in the circadian system in Scenario 3. (a-c) The change of the RMSE and the NCR upon the noise magnitude in Scenario 3. In (a), only the noise in \textit{Process P} was considered. In (b), only the noise in \textit{Process L} was considered. In (c), the two noises were considered together.}
    \label{fig:figure4.2.}
\end{figure}

The circadian clock state $\mathbf{x}$ is influenced by noise from various sources \cite{takahashi2017transcriptional}. For instance, stochastic biochemical reactions occurring in the circadian clock and many biological systems interacting with the circadian clock, such as a metabolic system, cause the randomness of the clock state. This can be described with the process noise matrix $\mathbf{K}$ \cref{eq::dXt} of the form 
\begin{equation}
    \mathbf{K} = 
    \begin{pmatrix} 
    \sigma_{P}^2 & 0 & 0 \\ 
    0 & \sigma_{P}^2 & 0 \\ 
    0 & 0 & \sigma_{L}^2 
    \end{pmatrix}
    \label{eq::system-noise}
\end{equation}
where $\sigma_{P}$ represents the magnitude of the noise directly affecting the clock (\textit{Process P} in Figure \ref{fig::3.1}a) and $\sigma_{L}$ represents the magnitude of the noise in the biological process (\textit{Process L} in Figure \ref{fig::3.1}a) that transmits the external light signal from the retina to the clock via the retinohypothalamic tract.\\
\indent Here, we studied the relationship between the process noise in the circadian clock and the phase estimate using our filtering algorithm. We first generated $in$-$silico$ data mimicking the typical human lifestyle illustrated in Figure \ref{fig::3.3}b. We set the sleep offset time $t_{i}^{\text{off}}$ and onset time $t_{i}^{\text{on}}$ of day $i$ as $7$am and $11$pm, respectively, so that activity level on day $i$ $a_{i}(t)=0$ if $t\in \left[0, \hspace{0.12cm} t_{i}^{\text{off}}\right) \cup \left[t_{i}^{\text{on}}, \hspace{0.12cm} 24\right]$. We next defined $a_{i}(t)$ in wake time $\left[t_{i}^{\text{off}}, \hspace{0.12cm} t_{i}^{\text{on}}\right)$ as done in previous work \cite{dijk2012amplitude, mott2011model}. Specifically, we subdivided the wake time into three stages: morning $\left[t_{i}^{\text{off}}, \hspace{0.12cm} t_{i}^{\text{off}}+5\right)$, afternoon $\left[t_{i}^{\text{off}}+5, \hspace{0.12cm} t_{i}^{\text{off}}+10\right)$, and evening $\left[t_{i}^{\text{off}}+10, \hspace{0.12cm} t_{i}^{\text{on}}\right)$. Then, under the assumption that level of external signals (e.g., light or activity) is low in the morning, high in the afternoon, and again low in the evening like the change of light intensity during the course of a day, $a_{i}(t)$ was set as 
\begin{equation}
    a_{i}(t) = \begin{cases} 
    \max(\mu_{l}+N(0,\sigma_{l}^2),0) & \text{if } t_{i}^{\text{off}} \leq t < t_{i}^{\text{off}}+5 \\
    \max(\mu_{h}+N(0,\sigma_{h}^2),0) & \text{if } t_{i}^{\text{off}}+5 \leq t < t_{i}^{\text{off}}+10 \\
    \max(\mu_{l}+N(0,\sigma_{l}^2),0) & \text{if } t_{i}^{\text{off}}+10 \leq t < t_{i}^{\text{on}}  \\
    0 & \text{otherwise}
    \end{cases}
\end{equation}
where $\mu_{l}$ and $\sigma_{l}$ denote the mean and the standard deviation of the low activity level, respectively, and $\mu_{h}$ and $\sigma_{h}$ denote the mean and the standard deviation of the high activity level, respectively. We took $a_{i}(t)$ as the maximum between a sample from $\mu + N(0,\sigma^2)$ and 0 to prevent it from being negative. Note that we set $\mu_{l}$ and $\mu_{h}$ as $5$ steps/min and $25$ steps/min, respectively, with which the activity signal is converted on average to a typical ordinary (500 lux) - bright (2000 lux) - ordinary (500 lux) light exposure by \cref{eq::acttolux}. Under this activity setting, we simulated the HR rhythm with a nadir at the midpoint of sleep $3$am using the model \cref{eq::HR} based on the previous study \cite{bowman2021method}. Then, we finally defined the internal phase $\phi_i^{true}$ to be $4$am based on its relationship with the HR phase \cref{eq::HRphase}. More details of the simulated setting named Scenario 1 are given in Table \ref{tab::summary_scenario_case}.\\
\indent We applied our method to Scenario 1 and analyzed the estimates with the two quantities, root mean square error (RMSE) and non-coverage rate (NCR). 
\begin{table}[hbt!]
    \centering
    \begin{tabular}{
        l
        S[table-format = 3]
        S[table-format = 2]
        S[table-format = 1.3]
        S[table-format = -2.2]
        S[table-format = 1.3]
        S[table-format = 1.3]
        S[table-format = 2.2]
        }
        \toprule
        \multicolumn{3}{r}{Summary of Scenarios}\\
        \cmidrule(lr){2-4}        
         {Parameters} & {Scenario 1} & {Scenario 2} & {Scenario 3}  \\
        \midrule
        {\textit{Activity level}} \\
        {morning $\&$ evening $\mu_{l}$} & {5 steps/min} & {5 steps/min} & {5 steps/min}\\ 
        {afternoon $\mu_{h}$} & {25 steps/min} & {25 steps/min} & {25 steps/min}\\
        {sleep} & {0 steps/min} & {0 steps/min} & {0 steps/min} \\ \\
        {\textit{Activity uncertainty}} \\
        {$\sigma_l$} & {7.5 steps/min} & {7.5 steps/min} & {7.5 steps/min} \\
        {$\sigma_h$} & {30 steps/min} & {30 steps/min} & {30 steps/min} \\
        {$\sigma_s$} & {0 steps/min} & {0 steps/min} & {2 steps/min} \\
        {$\sigma_t$} & {0 hrs} & {1.5 hrs} & {1.5 hrs} \\
        {HR phase} & {3hr} & {3hr} & {3hr} \\
        {Internal phase} & {4hr} & {4hr} & {4hr} \\ \\ 
        {\textit{HR signal}} \\
        {$\mu$} & {70 bpm} & {70 bpm} & {70 bpm} \\
        {$a$} & {4 bpm} & {4 bpm} & {4 bpm} \\
        {$\phi^{HR}$} & {3hr} & {3hr} & {3hr} \\
        {$d$} & {0.3} & {0.3} & {0.3} \\ \\
        {\textit{HR uncertainty}} \\
        {$\sigma_\epsilon$} & {3 bpm} & {3 pm} & {7 bpm} \\
        {$\alpha$} & {0.9} & {0.9} & {0.95} \\
        \bottomrule \\
    \end{tabular}
    \caption{Summary of the simulated scenarios. In Scenario 1, the randomness of activity level in wake time is only considered. In Scenario 2, the randomness of the sleep onset and offset times is also considered. In Scenario 3, the small randomness of activity level in sleep time, due to tossing and turning or measurement noise of wearables, for example, is additionally considered. Moreover, the HR uncertainty is increased to account for the potentially large measurement errors in wearable HR data in real-world settings.}
    \label{tab::summary_scenario_case}
\end{table}
The RMSE is the standard deviation of the difference between the mean phase estimate $\bar{\phi}_{i|i}$ and the true phase $\phi^{true}_{i}$, which is defined to be 
\begin{equation}
    RMSE = \sqrt{\frac{1}{N} \sum_{i=1}^{N} e_i^2}
    \label{eq::RMSE}
\end{equation}
where $e_i=\left|\mathbb{E}[\phi_{i|i}] - \phi_i^{true}\right|$ and $N$ is the total number of days that the phase was estimated. The NCR is the rate of days that the true phase is not included in the 95$\%$ credible interval, which is defined to be
\begin{equation}
    NCR = 1- \frac{1}{N} \sum_{i=1}^{N} \mathbbm{1}^{i}_{\text{CI}}(\phi_{i}^{true})
    \label{eq::ncr}
\end{equation}
where $\mathbbm{1}^{i}_{\text{CI}}$ denotes the indicator function that maps inputs to one if they are included in the 95$\%$ credible interval on day $i$. Note that our credible estimates are based on a Bayesian procedure, and hence we do not focus on the agreement between the nominal and empirical coverages. Instead, we just focus on whether the true phase on each day lies within the corresponding credible interval estimate. The graphical illustration of the RMSE and the NCR is shown in Figure \ref{fig::4.1}a. The small value of both the RMSE and the NCR means that both the mean and the uncertainty of the phase estimate are accurately identified. This interpretation of the relationship between the two quantities and the phase estimate is shown in Figure \ref{fig::4.1}b.\\
\indent Using these two measures, we analyzed the phase estimates obtained by applying our method to Scenario 1. Specifically, we first let $\sigma_{L}=0$ and took $\mathbf{K}$ of the form \cref{eq::p-noise} to focus on exploring the relationship between the phase estimate and the noise directly affecting the pacemaker (i.e., the noise in \textit{Process P}). Then, we varied
\begin{equation}
    \mathbf{K} = 
    \begin{pmatrix} 
    \sigma_{P}^2 & 0 & 0 \\ 
    0 & \sigma_{P}^2 & 0 \\ 
    0 & 0 & 0 
    \end{pmatrix}.
    \label{eq::p-noise}
\end{equation}
$\sigma_{P}$ of \cref{eq::p-noise} and calculated the RMSE and the NCR as shown in Figure \ref{fig::4.1}c. In Figure \ref{fig::4.1}c, small $\sigma_{P}$ results in the small RMSE but the large NCR. This indicates that ruling out the noise in the pacemaker might cause problems related to estimating the phase uncertainty. On the other hand, when $\sigma_{P}$ is set to be large, the NCR is small, but the RMSE is large, showing that overrating $\sigma_{P}$ can lead to an inaccurate estimation of the mean phase. Accordingly, $\sigma_{P}$ of appropriate magnitude is required for accurate estimation of both the mean and the uncertainty of the phase. We similarly analyzed the relationship between the phase estimate and the noise in the light processor (i.e., the noise in \textit{Process L}) with $\textbf{K}$ of the form
\begin{equation}
    \mathbf{K} = 
    \begin{pmatrix} 
    0 & 0 & 0 \\ 
    0 & 0 & 0 \\ 
    0 & 0 & \sigma_{L}^2 
    \end{pmatrix}.
    \label{eq::l-noise}
\end{equation}
Although $\sigma_{L}$ changes, the RMSE and the NCR remain small and large, respectively, as shown in Figure \ref{fig::4.1}d, again indicating that the noise in \textit{Process P} needs to be considered for accurate phase estimation. Moreover, the weak dependence of the RMSE and the NCR on $\sigma_{L}$ suggests that the noise in \textit{Process L} affects the state of the molecular clocks more weakly than that in \textit{Process P}. This is supported by the fact that the noise in \textit{Process L} originated from the light signal transduction pathway (not the molecular clocks), so it indirectly affects the clock state. To study this, we explored the association of the phase estimate with the total noise in \textit{Process P} and \textit{Process L} by taking $\mathbf{K}$ of the form 
\begin{equation}
    \mathbf{K} = \sigma_{\mathbf{K}}^2 \cdot I_3
\end{equation}
where $\sigma_{\mathbf{K}}$ denotes the magnitude of the total noise, and $I_3$ is the identity matrix of order $3$. As expected, the relationship between the phase estimate and the total process noise shown in \ref{fig::4.1}e is mainly determined by the influence of the noise in \textit{Process P} shown in Figure \ref{fig::4.1}c. These results obtained using our method provide an understanding of the relationship between the molecular clock state and noise, which is experimentally unobservable. Importantly, based on this knowledge, we can infer the magnitude of the process noise that accurately describes the mean and the uncertainty of the circadian phase, leading to both the small RMSE and NCR. For example, from Figure \ref{fig::4.1}e, we can figure out that $\sigma_{\mathbf{K}}$ yielding both the small RMSE (0.164hr) and NCR (0.047) is $10^{-3}$ in Scenario 1.\\
\indent We next applied our method to more realistic $in$-$silico$ data named Scenario 2. In Scenario 2, the settings of Scenario 1 are replicated except that uncertainty is introduced into sleep offset time $t_i^{\text{off}}$ and onset time $t_i^{\text{on}}$ as follows:
\begin{equation}
    t_i^{\text{off}} \sim 7+N(0,\sigma_{t}^2), \quad t_i^{\text{on}} \sim 23+N(0,\sigma_{t}^2)
    \label{eq::t_onset_offset}
\end{equation}
where $\sigma_{t}=1.5$hr represents the uncertainty of the sleep onset and offset times (see Table \ref{tab::summary_scenario_case} for more details). Figures \ref{fig::4.1}f, \ref{fig::4.1}g, and \ref{fig::4.1}h describe the association of the phase estimate with the noise in \textit{Process P}, the noise in \textit{Process L}, and the total noise, respectively, in Scenario 2, like Figures \ref{fig::4.1}c, \ref{fig::4.1}d, and \ref{fig::4.1}e. In Scenario 2, the relationship between phase estimate and noise is mainly governed by noise in \textit{Process P}, as in Scenario 1. However, unlike Scenario 1, the change of the RMSE upon the addition of $\sigma_P$ is not large, indicating the discrepancy of the circadian dynamics in different sleep/wake patterns.\\
\indent We finally analyzed the most challenging but realistic $in$-$silico$ data named Scenario 3. In Scenario 3, the settings of Scenario 2 are repeated except that the randomness of activity level in sleep time, which originated from small unconscious movements during sleep and measurement noise of wearables, is introduced as follows:
\begin{equation}
    a_{i}(t) = \max(N(0,\sigma_{s}^2),0) \hspace{0.4cm} \text{if } 0\leq t < t_{i}^{\text{off}} \hspace{0.2cm} \text{or} \hspace{0.2cm} t_{i}^{\text{on}} \leq t \leq 24 \\
\end{equation}
where $\sigma_{s}=2$ steps/min represents the magnitude of the small randomness of activity level in sleep time. Moreover, we increased the uncertainty of HR measurements to account for the potentially large measurement noise in HR data collected in real-world settings reported on \cite{bowman2021method}. See Table \ref{tab::summary_scenario_case} for more details. Figures \ref{fig::4.2}a, \ref{fig::4.2}b, and \ref{fig::4.2}c present the association of the phase estimate with the noise in \textit{Process P}, the noise in \textit{Process L}, and the total noise, respectively, in Scenario 3, like Figure \ref{fig::4.1}. As in Scenario 1 and 2, the relationship is mainly determined by the noise in \textit{Process P}. However, unlike the other scenarios, the RMSE and the NCR are very large even if $\sigma_{\mathbf{P}}$ is small, as shown in Figure \ref{fig::4.2}a. Moreover, unlike Scenario 1, the RMSE is large when only the noise in \textit{Process L} is considered as shown in Figure \ref{fig::4.2}b. This indicates that the noise in \textit{Process P} needs to be necessarily considered in realistic settings to estimate both the mean and uncertainty of the phase accurately. Importantly, even in this challenging scenario, our method can estimate the phase within $1$hr (i.e., $\text{RMSE}<1$hr) with a carefully chosen $\sigma_{\mathbf{K}}=6\cdot 10^{-3}$ as highlighted by an arrow in Figure \ref{fig::4.2}c, which is impossible with the previous methods \cite{huang2021predicting, bowman2021method} (see section \ref{sec::4.2} below for more details).
\subsection{Performance comparison} 
\label{sec::4.2}
\begin{figure}[hbt!]
    \label{fig::4.3}
    \centering
    \includegraphics[scale = 0.49]{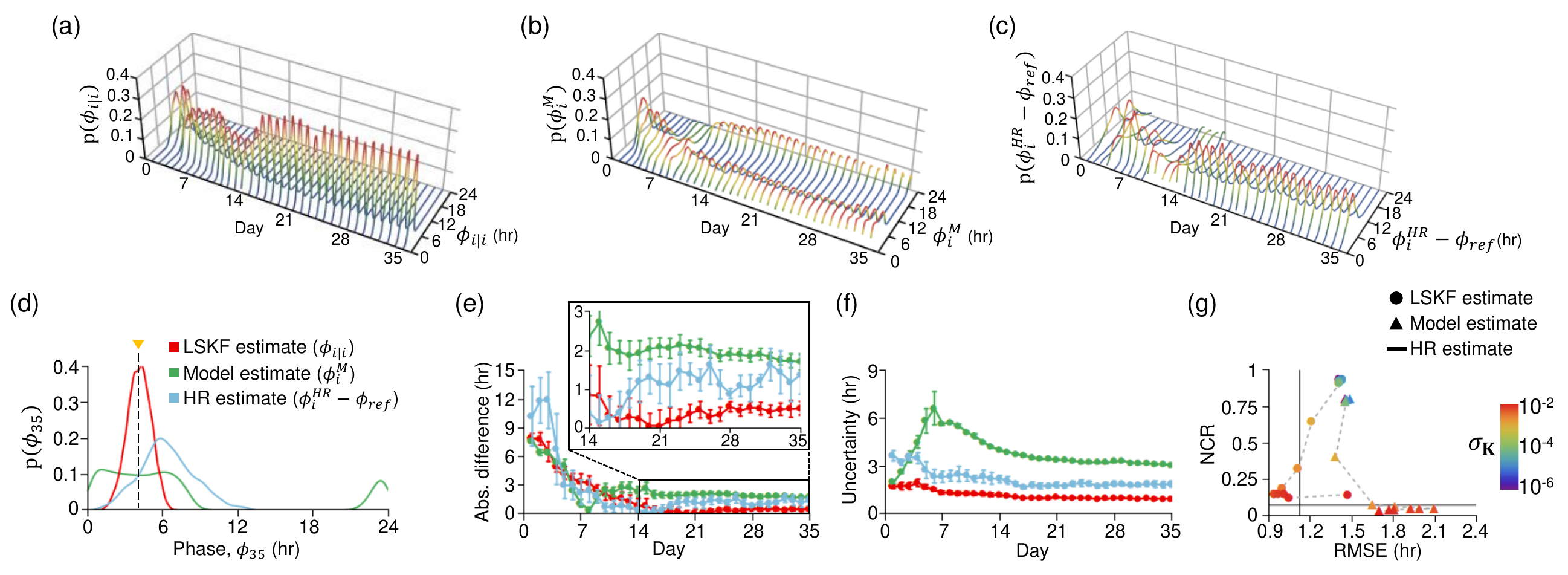}
    \caption{Our data assimilation approach can improve the accuracy of the circadian phase estimation from wearable measurements. (a-c) The evolution of the posterior distributions of the phase that are obtained using our LSKF method (a), only using the mathematical model \cref{eq::dXt} \cite{huang2021predicting} (b), and only using the HR phase estimate \cite{bowman2021method} (c). Here, the three methods were applied to Scenario 3 and $\sigma_{\mathbf{K}}$ of $\mathbf{K}=\sigma_{\mathbf{K}}^2\cdot I_3$ was set as $6\cdot 10^{-3}$ because it allows for accurately capturing the circadian dynamics in Scenario 3 as shown in Figure \ref{fig::4.2}c. (d) The posterior distributions estimated by the three methods on the last day. (e) The absolute difference between the true phase and the mean phase estimate over days. Here, five independent $in$-$silico$ data of Scenario 3 were analyzed, and the five computed absolute differences were averaged. Error bars denote the standard error of the mean (SEM). (f) The phase uncertainty (i.e., the standard deviation of the phase estimate) over days. The mean and SEM of the phase uncertainties on each day were calculated as in (e). (g) The RMSE and the NCR of the three methods are computed with various $\sigma_{\mathbf{K}}$. Note that the estimation result of the method that computes the clock phase by subtracting $\phi_{ref}$ from the HR phase estimate \cite{bowman2021method} is independent of $\sigma_{\mathbf{K}}$. Thus, the RMSE and the NCR of the method are represented as the vertical and horizontal lines, respectively.}
    \label{fig:figure4.3.}
\end{figure}

Here, we compared the performance of our LSKF method with the recently developed methods that do not use data assimilation to demonstrate its capabilities for tracking the circadian phase from wearable data. Specifically, two previous methods were adopted for the comparison. The first previous method proposed in \cite{huang2021predicting} estimates the clock state based solely on the mathematical model of the human circadian clock \cref{eq::dXt} that takes activity wearable data as an input as described in section \ref{sec::3.3.2}. The second one first estimates the HR phase $\phi_{i}^{HR}$ by fitting \cref{eq::HR} to wearable data using Goodman and Weare's affine-invariance MCMC method and then compute the clock state by subtracting $\phi_{ref}=-1$ from $\phi_{i}^{HR}$ \cite{bowman2021method}. These two methods were chosen for the comparison study for the following reasons: (i) they are recently developed and thus sophisticated; (ii) our LSKF method integrates the prediction from the model used in \cite{huang2021predicting} with the physiological parameter (i.e., the HR phase) obtained using the method in \cite{bowman2021method}. We applied the methods to the $in$-$silico$ data of Scenario 3, which is the most realistic but challenging setting and calculated the absolute error (i.e., $e_i$ in \cref{eq::RMSE}) and the standard deviation (i.e., the uncertainty) of the estimate on each day. Then, we compared the estimates of the two previous methods \cite{bowman2021method, huang2021predicting} with those computed with our LSKF method. Figures \ref{fig::4.3}a, \ref{fig::4.3}b, and \ref{fig::4.3}c show the evolution of the posterior distributions of the phase over days estimated using our LSKF method, only using the mathematical model \cref{eq::dXt} \cite{huang2021predicting}, and only using the HR phase estimate \cite{bowman2021method}, respectively. They show that the distributions estimated using our method are narrower than those estimated using the previous methods. Moreover, the filtering algorithm is more accurate than the previous methods, as shown in Figures \ref{fig::4.3}d. Indeed, the absolute error and the uncertainty of our method are smaller than those of the others, as shown in Figures \ref{fig::4.3}e and \ref{fig::4.3}f. Our method can estimate the phase within $1$hr ($\text{RMSE}=0.942$hr) while the others cannot (model estimate \cite{huang2021predicting}: $\text{RMSE}=1.996$hr and HR estimate \cite{bowman2021method}: $\text{RMSE}=1.527$hr). Notably, the improvement is preserved although another widely used mathematical model of the human circadian clock \cite{jewett1999revised} is adopted instead of the FJK model \cite{forger1999simpler}, as described in Supplementary Materials. This demonstrates the value of our data assimilation technique.\\
\indent We next investigated whether the outperformance of our method is preserved even if the magnitude of the process noise $\sigma_{\mathbf{K}}$ changes. Specifically, in Figures \ref{fig::4.3}e and \ref{fig::4.3}f, $\sigma_{\mathbf{K}}=6\cdot 10^{-3}$ was used because it allows our method to accurately capture the dynamics of the circadian systems, resulting in the phase estimation with the small RMSE and NCR in Scenario 3 as shown in Figure \ref{fig::4.2}c. Importantly, even if this $\sigma_{\mathbf{K}}$ varies, and thus the performance of our method becomes worse, its RMSE is overall smaller than the RMSE of the previous methods, as presented in Figure \ref{fig::4.3}g. In particular, the RMSE of our method is always smaller than that of the method based solely on the model prediction. Moreover, the NCR of our method is comparable with that of the previous methods with a carefully selected $\sigma_{K}$. These results demonstrate the benefits of the filtering approach on the circadian phase estimation from wearable data.\\
\indent Lastly, using our method, we analyzed scenarios in which individuals experience misalignment between activity-rest rhythms and external light-dark cycles (Supplementary Materials). Indeed, the estimation results do not change significantly depending on the model input (activity vs light), as shown in Supplementary Materials, supporting the use of activity measurements when light information is not available.
\section{Application to a real data set}
\label{sec::real data}
\begin{figure}[hbt!]
    \label{fig::5.1}
    \centering
    \includegraphics[scale = 0.74]{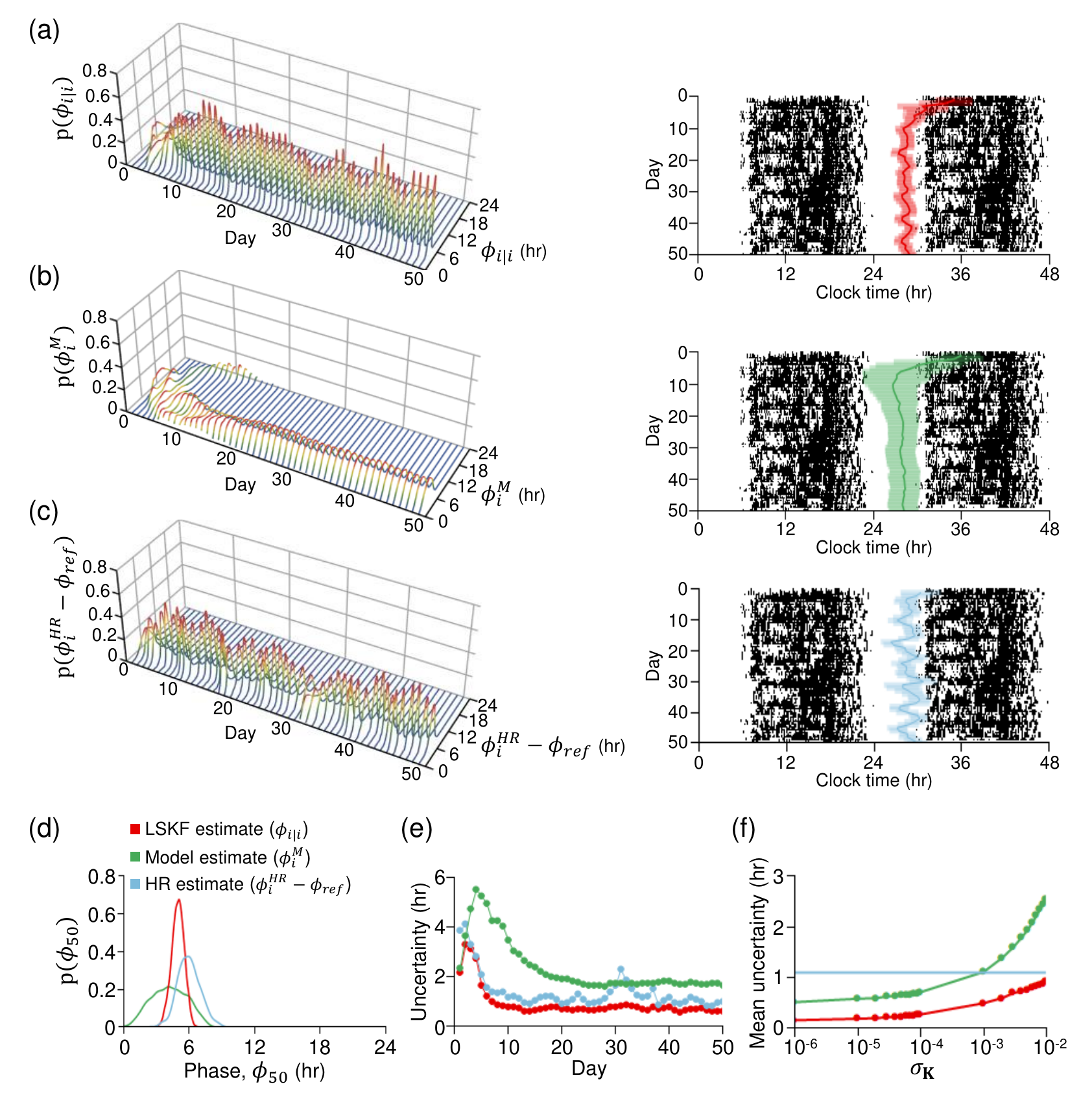}
    \caption{Estimation of the circadian phase from real-world wearable-device data. (a-c) The phase estimates obtained using our filtering method (a), based solely on model prediction (b), or by direct calculation from the HR phase (c). In the left panel, the evolution of the posterior distributions of the phase is presented. In the right panel, their mean and standard deviation are plotted together with actograms. Here, $\sigma_{\mathbf{K}}$ of $\mathbf{K}=\sigma_{\mathbf{K}}\cdot I_3$ was set as $6\cdot 10^{-3}$ that allows the accurate identification of the circadian dynamics in Scenario 3 (Figure \ref{fig::4.3}c). (d) The posterior distributions estimated by the three methods on the last day. (e) The phase uncertainty over days. (f) The average of the phase uncertainties over days that is computed with various $\sigma_{\mathbf{K}}$. Note that the estimation result of the method that directly calculates the molecular clock state only from the HR phase is independent of $\sigma_{\mathbf{K}}$. Thus, the uncertainty computed using the method is constant over $\sigma_{\mathbf{K}}$.}
    \label{fig:figure5.1.}
\end{figure}

We developed a publicly available computational package \href {https://github.com/phillee62/LSKF_circadian} {(https://github.com/phillee62/LSKF$\_$circadian)} to facilitate the use of our algorithm. Using this computational package, we analyzed the previously collected real-world wearable HR and activity data \cite{bowman2021method} to further demonstrate the value of our algorithm in Figure \ref{fig::5.1}. Specifically, we estimated the circadian phase using our method with the value of $\sigma_{\mathbf{K}}=6\cdot 10^{-3}$ that accurately captures the circadian dynamics in the most challenging but realistic simulation scenario, Scenario 3, as shown in Figure \ref{fig::4.2}c. Then, the estimated phases were compared with the phases that were estimated based solely on the mathematical model taking in inputs of activity measurements as the previous studies did \cite{huang2021predicting}. We also compared our estimates with those calculated by subtracting $\phi_{ref}=-1$ from the HR phase estimate ($\phi_{i}^{HR}$ in \cref{eq::HRphase}), which was obtained by fitting \cref{eq::HR} to wearable HR and activity data using Goodman and Weare's affine-invariance MCMC method \cite{bowman2021method}. Figures \ref{fig::5.1}a, \ref{fig::5.1}b, and \ref{fig::5.1}c show the phases estimated by our method, the predicted phases by the mathematical model, and the phases directly calculated from the HR phases, respectively. In their left panels, the evolution of the posterior distributions of the phase is presented. The distributions on the last day are presented in Figure \ref{fig::5.1}d. In the right panels of Figures \ref{fig::5.1}a, \ref{fig::5.1}b, and \ref{fig::5.1}c, the mean and the standard deviation of the distributions are shown with double-plotted actograms. The standard deviation representing the uncertainty of the phase estimate is quantified in Figure \ref{fig::5.1}e. These figures show that the posterior distributions computed using our LSKF method are narrower than the others. Figure \ref{fig::5.1}f shows that this pattern is preserved even if the magnitude of the process noise $\sigma_{\mathbf{K}}$ is varied so that the noise magnitude leading to the best performance of our method is not exploited in the estimation. This indicates that the uncertainty of the phase estimate can be reduced when combining the model prediction and the measurements from wearables, which is consistent with our results in the numerical experiments. Moreover, the phase estimated by our method differs from the others, for example, as shown in  Figure \ref{fig::5.1}d. Considering this with the outperformance of our method in the numerical study, filtering approaches like our algorithm might be needed for accurate circadian phase estimation.
\section{Conclusions}
\label{sec::conclusion}
We developed a Kalman filter framework for estimating the state space of the molecular clocks in tissues. As well as estimating the state space, our method can quantify its uncertainty systematically, which is impossible with the previous methods based on ODE models \cite{huang2021predicting, stone2019application, woelders2017daily}. Numerical experiments in Figure \ref{fig::4.3} showed how much the uncertainty could be reduced when utilizing the indirect observation of the molecular clock state. Moreover, our method overall outperforms the previous methods, \cite{bowman2021method, huang2021predicting, mayer2022consumer} as presented in Figure \ref{fig::4.3}. These results suggest new avenues for exploiting noninvasive wearable data for chronotherapy in free-living conditions.\\
\indent In addition to its accuracy, our framework can successfully minimize the influence of initial conditions on the phase estimation. Our method continuously updates the phase estimate by the measurement-update step described in section \ref{sec::3.4.2}. This results in the rapid identification of the true phase even if arbitrary initial conditions are given, as shown in Figure \ref{fig::3.3}c. Specifically, when a large initial covariance matrix to account for arbitrary initial conditions is given, cubature gain, denoted as $\mathbf{W}$ in Algorithm \ref{alg::measurement-update}, becomes large. This leads to an adjustment of the phase estimate toward the measurement (i.e., the HR phase), which is independent of the initial conditions of the clock state. This algorithmic procedure rapidly negates the undesired influence of initial conditions, which is not available in the previous methods solely relying on the convergence nature of a van der Pol limit cycle \cite{huang2021predicting, walch2019sleep}. Importantly, in the analysis of real data, our method more rapidly negated the influence of initial conditions ($<$ 10 days) than the method based solely on the model prediction, as shown in Figures \ref{fig::5.1}a, \ref{fig::5.1}b, and \ref{fig::5.1}e. This indicates the applicability of our method in real-life settings, where minimizing the effect of the initial conditions is desired.\\
\indent In this study, a recently developed new state space estimation method called the level set Kalman filter \cite{wang2021level} was used to integrate the model prediction with the physiological proxies extracted from wearable data. To the best of our knowledge, our work is the first attempt to use a Kalman filter to estimate the state space of the molecular clocks in tissues in real-life settings from wearable data. It has been previously used to reduce the computational cost of performing the maximum likelihood estimation of physiological parameters of core body temperature rhythms \cite{brown19949}. Specifically, the standard discrete-time Kalman filter was exploited to efficiently compute the Cholesky factor of the inverse of the covariance matrix of maximum likelihood estimates. Considering this, our study first suggests the potential of data assimilation approaches based on the Kalman filter for accurate circadian phase estimation from wearable data. An interesting avenue for future research is to identify suitable filtering methods for the phase estimation like the LSKF. Another important future work is to test ability of different types of steps-to-light function, such as a continuous conversion function, to estimate the circadian phase.\\
\indent Our filtering problem was tackled by using the recently proposed LSKF \cite{wang2021level}. It tracks a Gaussian level set over time to solve the local linear approximation of the  Fokker-Planck equation instead of calculating the moments of distribution as shown in Figure \ref{fig::3.2}a. This novel approach to filtering problems can be applied to solve the chemical Langevin equation, a diffusion approximation of the chemical master equation describing stochastic dynamics of chemical systems \cite{schnoerr2017approximation}. It would be interesting in future work to compare its performance with existing methods such as the system size expansion \cite{van1992stochastic} and moment closure approximations \cite{gillespie2009moment}.\\
\indent Our LSKF-based filtering approach also has broad applicability in biological studies. Specifically, it can  estimate the experimentally unobservable (i.e., hidden) biochemical parameters \cite{forger2017biological, kim2022systematic, sun2008extended}. For instance, the estimation problem of unobservable time-varying protein production rate $\mathbf{m}^{1}_{t}$ and degradation rate $\mathbf{m}^{2}_{t}$ can be formulated as a filtering problem with the model of the form 
\begin{equation}
\begin{aligned}
    d \mathbf{m}_t &= \mathbf{v}(\mathbf{m}_t, t) dt + \sqrt{\mathbf{K}} dW_t\\
    \dot{\mathbf{p}}_k &= \mathbf{m}^{1}_{t_k} - \mathbf{m}^{2}_{t_k}\mathbf{p}_k  + \epsilon_k
    \label{eq::bioe}
\end{aligned}
\end{equation}
where $\mathbf{m}_{t} = \begin{bmatrix} \mathbf{m}^{1}_{t} \hspace{0.1cm} \mathbf{m}^{2}_{t} \end{bmatrix}^T$ whose dynamics is described with a nonlinear drift function $\mathbf{v}$, $\mathbf{p}_k$ and $\dot{\mathbf{p}}_k$ denote the protein abundance and its derivative, respectively, that can be measured with experimental techniques such as bioluminescence \cite{paley2014bioluminescence}. Note that $\mathbf{v}$ is chosen depending on the protein of interest.\\
\indent Our method can estimate the circadian phase with low uncertainty from wearables compared with the others as shown in Figures \ref{fig::4.3}f and \ref{fig::5.1}e. This benefit is preserved for a range of the parameter $\sigma_{\mathbf{K}}$ of the process noise matrix, as shown in Figure \ref{fig::5.1}f. This demonstrates the capabilities of our method for tracking the circadian phase. In the analysis of the real-world data (Figure \ref{fig::5.1}), we used the process noise matrix $\mathbf{K} = \sigma_{\mathbf{K}}\cdot I_3$ with which our algorithm can accurately estimate the circadian phase from the most challenging but realistic \textit{in-silico} data (Figure \ref{fig::4.2}c). This approach based on \textit{in-silico} tests has limitations because there is accumulating evidence for a large inter-and intraindividual variability in circadian variables \cite{kim2020wearable}. Thus, future work is needed to develop a more systematic and rigorous approach that tailors $\mathbf{K}$ and the other model parameters to individual circadian physiology. One promising approach is to use methods, such as the output correlation approach, based on relations between the noise parameters and the covariance function of observable measurements \cite{bulut2011process}.\\
\indent Accurate estimation of circadian phase and its uncertainty enables personalized and real-time monitoring of progression of various diseases, including viral, bacterial, and neurodegenerative diseases \cite{li2018fractal, li2017digital, mayer2022consumer, natarajan2020assessment}. For instance, it has been recently shown that COVID-19 symptom onset correlates with an increased circadian phase uncertainty in \cite{mayer2022consumer}. Thus, the phase and its uncertainty accurately estimated by our method can be exploited for early disease diagnosis in real-world settings. Furthermore, accurate circadian phase estimation benefits sleep scoring based on wearable data. Specifically, it has been shown that the phase extracted only using the mathematical model can increase the accuracy of the sleep scoring machine-learning algorithm in \cite{walch2019sleep}. Thus, more abundant and accurate information about the circadian phase is expected to enhance the accuracy of the existing wearable-based sleep scoring algorithms. This promising clinical applicability demonstrates that our filtering approach may provide an important advance in precision medicine in real-life conditions.

\bibliographystyle{siamplain}
\bibliography{references}

\begin{thebibliography}{10}

\bibitem{arasaratnam2009cubature}
{\sc I.~Arasaratnam and S.~Haykin}, {\em Cubature kalman filters}, IEEE
  Transactions on automatic control, 54 (2009), pp.~1254--1269.

\bibitem{arasaratnam2010cubature}
{\sc I.~Arasaratnam, S.~Haykin, and T.~R. Hurd}, {\em Cubature kalman filtering
  for continuous-discrete systems: theory and simulations}, IEEE Transactions
  on Signal Processing, 58 (2010), pp.~4977--4993.

\bibitem{becker2019time}
{\sc L.~Becker and N.~Rohleder}, {\em Time course of the physiological stress
  response to an acute stressor and its associations with the primacy and
  recency effect of the serial position curve}, PLoS One, 14 (2019),
  p.~e0213883.

\bibitem{bowman2021method}
{\sc C.~Bowman, Y.~Huang, O.~J. Walch, Y.~Fang, E.~Frank, J.~Tyler, C.~Mayer,
  C.~Stockbridge, C.~Goldstein, S.~Sen, et~al.}, {\em A method for
  characterizing daily physiology from widely used wearables}, Cell reports
  methods, 1 (2021), p.~100058.

\bibitem{bressloff2014stochastic}
{\sc P.~C. Bressloff}, {\em Stochastic processes in cell biology}, vol.~41,
  Springer, 2014.

\bibitem{brown2000statistical}
{\sc E.~N. Brown, Y.~Choe, H.~Luithardt, and C.~A. Czeisler}, {\em A
  statistical model of the human core-temperature circadian rhythm}, American
  Journal of Physiology-Endocrinology and Metabolism, 279 (2000),
  pp.~E669--E683.

\bibitem{brown1992statistical}
{\sc E.~N. Brown and C.~A. Czeisler}, {\em The statistical analysis of
  circadian phase and amplitude in constant-routine core-temperature data},
  Journal of Biological Rhythms, 7 (1992), pp.~177--202.

\bibitem{brown1999statistical}
{\sc E.~N. Brown and H.~Luithardt}, {\em Statistical model building and model
  criticism for human circadian data}, Journal of Biological Rhythms, 14
  (1999), pp.~609--616.

\bibitem{brown19949}
{\sc E.~N. Brown and C.~H. Schmid}, {\em Application of the kalman filter to
  computational problems in statistics}, in Methods in enzymology, vol.~240,
  Elsevier, 1994, pp.~171--181.

\bibitem{buckert2014acute}
{\sc M.~Buckert, C.~Schwieren, B.~M. Kudielka, and C.~J. Fiebach}, {\em Acute
  stress affects risk taking but not ambiguity aversion}, Frontiers in
  neuroscience, 8 (2014), p.~82.

\bibitem{bulut2011process}
{\sc Y.~Bulut, D.~Vines-Cavanaugh, and D.~Bernal}, {\em Process and measurement
  noise estimation for kalman filtering}, in Structural Dynamics, Volume 3:
  Proceedings of the 28th IMAC, A Conference on Structural Dynamics, 2010,
  Springer, 2011, pp.~375--386.

\bibitem{dibner2010mammalian}
{\sc C.~Dibner, U.~Schibler, and U.~Albrecht}, {\em The mammalian circadian
  timing system: organization and coordination of central and peripheral
  clocks}, Annual review of physiology, 72 (2010), pp.~517--549.

\bibitem{diekman2018reentrainment}
{\sc C.~O. Diekman and A.~Bose}, {\em Reentrainment of the circadian pacemaker
  during jet lag: East-west asymmetry and the effects of north-south travel},
  Journal of theoretical biology, 437 (2018), pp.~261--285.

\bibitem{dijk2012amplitude}
{\sc D.-J. Dijk, J.~F. Duffy, E.~J. Silva, T.~L. Shanahan, D.~B. Boivin, and
  C.~A. Czeisler}, {\em Amplitude reduction and phase shifts of melatonin,
  cortisol and other circadian rhythms after a gradual advance of sleep and
  light exposure in humans}, PloS one, 7 (2012), p.~e30037.

\bibitem{forger2017biological}
{\sc D.~B. Forger}, {\em Biological clocks, rhythms, and oscillations: the
  theory of biological timekeeping},  (2017).

\bibitem{forger1999simpler}
{\sc D.~B. Forger, M.~E. Jewett, and R.~E. Kronauer}, {\em A simpler model of
  the human circadian pacemaker}, Journal of biological rhythms, 14 (1999),
  pp.~533--538.

\bibitem{forger2002reconciling}
{\sc D.~B. Forger and R.~E. Kronauer}, {\em Reconciling mathematical models of
  biological clocks by averaging on approximate manifolds}, SIAM Journal on
  Applied Mathematics, 62 (2002), pp.~1281--1296.

\bibitem{forger2004starting}
{\sc D.~B. Forger and D.~Paydarfar}, {\em Starting, stopping, and resetting
  biological oscillators: in search of optimum perturbations}, Journal of
  theoretical biology, 230 (2004), pp.~521--532.

\bibitem{forger2003detailed}
{\sc D.~B. Forger and C.~S. Peskin}, {\em A detailed predictive model of the
  mammalian circadian clock}, Proceedings of the National Academy of Sciences,
  100 (2003), pp.~14806--14811.

\bibitem{gillespie2009moment}
{\sc C.~S. Gillespie}, {\em Moment-closure approximations for mass-action
  models}, IET systems biology, 3 (2009), pp.~52--58.

\bibitem{goodman2010ensemble}
{\sc J.~Goodman and J.~Weare}, {\em Ensemble samplers with affine invariance},
  Communications in applied mathematics and computational science, 5 (2010),
  pp.~65--80.

\bibitem{hannay2019macroscopic}
{\sc K.~M. Hannay, V.~Booth, and D.~B. Forger}, {\em Macroscopic models for
  human circadian rhythms}, Journal of Biological Rhythms, 34 (2019),
  pp.~658--671.

\bibitem{hannay2018macroscopic}
{\sc K.~M. Hannay, D.~B. Forger, and V.~Booth}, {\em Macroscopic models for
  networks of coupled biological oscillators}, Science advances, 4 (2018),
  p.~e1701047.

\bibitem{hong2021personalized}
{\sc J.~Hong, S.~J. Choi, S.~H. Park, H.~Hong, V.~Booth, E.~Y. Joo, and J.~K.
  Kim}, {\em Personalized sleep-wake patterns aligned with circadian rhythm
  relieve daytime sleepiness}, Iscience, 24 (2021), p.~103129.

\bibitem{huang2021predicting}
{\sc Y.~Huang, C.~Mayer, P.~Cheng, A.~Siddula, H.~J. Burgess, C.~Drake,
  C.~Goldstein, O.~Walch, and D.~B. Forger}, {\em Predicting circadian phase
  across populations: a comparison of mathematical models and wearable
  devices}, Sleep, 44 (2021), p.~zsab126.

\bibitem{jewett1999revised}
{\sc M.~E. Jewett, D.~B. Forger, and R.~E. Kronauer}, {\em Revised limit cycle
  oscillator model of human circadian pacemaker}, Journal of biological
  rhythms, 14 (1999), pp.~493--500.

\bibitem{kalman1961new}
{\sc R.~E. Kalman and R.~S. Bucy}, {\em New results in linear filtering and
  prediction theory},  (1961).

\bibitem{kim2022systematic}
{\sc D.~W. Kim, H.~Hong, and J.~K. Kim}, {\em Systematic inference identifies a
  major source of heterogeneity in cell signaling dynamics: The rate-limiting
  step number}, Science advances, 8 (2022), p.~eabl4598.

\bibitem{kim2020wearable}
{\sc D.~W. Kim, E.~Zavala, and J.~K. Kim}, {\em Wearable technology and systems
  modeling for personalized chronotherapy}, Current Opinion in Systems Biology,
  21 (2020), pp.~9--15.

\bibitem{kronauer1990quantitative}
{\sc R.~Kronauer}, {\em A quantitative model for the effects of light on the
  amplitude and phase of the deep circadian pacemaker, based on human data},
  Sleep, 90 (1990), pp.~306--309.

\bibitem{kulikov2017accurate}
{\sc G.~Y. Kulikov and M.~V. Kulikova}, {\em Accurate continuous--discrete
  unscented kalman filtering for estimation of nonlinear continuous-time
  stochastic models in radar tracking}, Signal Processing, 139 (2017),
  pp.~25--35.

\bibitem{li2018fractal}
{\sc P.~Li, L.~Yu, A.~S. Lim, A.~S. Buchman, F.~A. Scheer, S.~A. Shea, J.~A.
  Schneider, D.~A. Bennett, and K.~Hu}, {\em Fractal regulation and incident
  alzheimer's disease in elderly individuals}, Alzheimer's \& Dementia, 14
  (2018), pp.~1114--1125.

\bibitem{li2017digital}
{\sc X.~Li, J.~Dunn, D.~Salins, G.~Zhou, W.~Zhou, S.~M.
  Sch{\"u}ssler-Fiorenza~Rose, D.~Perelman, E.~Colbert, R.~Runge, S.~Rego,
  et~al.}, {\em Digital health: tracking physiomes and activity using wearable
  biosensors reveals useful health-related information}, PLoS biology, 15
  (2017), p.~e2001402.

\bibitem{lovallo2006cortisol}
{\sc W.~R. Lovallo, N.~H. Farag, A.~S. Vincent, T.~L. Thomas, and M.~F.
  Wilson}, {\em Cortisol responses to mental stress, exercise, and meals
  following caffeine intake in men and women}, Pharmacology Biochemistry and
  Behavior, 83 (2006), pp.~441--447.

\bibitem{massin2000circadian}
{\sc M.~M. Massin, K.~Maeyns, N.~Withofs, F.~Ravet, and P.~G{\'e}rard}, {\em
  Circadian rhythm of heart rate and heart rate variability}, Archives of
  disease in childhood, 83 (2000), pp.~179--182.

\bibitem{mayer2022consumer}
{\sc C.~Mayer, J.~Tyler, Y.~Fang, C.~Flora, E.~Frank, M.~Tewari, S.~W. Choi,
  S.~Sen, and D.~B. Forger}, {\em Consumer-grade wearables identify changes in
  multiple physiological systems during covid-19 disease progression}, Cell
  Reports Medicine, 3 (2022), p.~100601.

\bibitem{mott2011model}
{\sc C.~Mott, G.~Dumont, D.~B. Boivin, and D.~Mollicone}, {\em Model-based
  human circadian phase estimation using a particle filter}, IEEE Transactions
  on Biomedical Engineering, 58 (2011), pp.~1325--1336.

\bibitem{natarajan2020assessment}
{\sc A.~Natarajan, H.-W. Su, and C.~Heneghan}, {\em Assessment of physiological
  signs associated with covid-19 measured using wearable devices}, NPJ digital
  medicine, 3 (2020), pp.~1--8.

\bibitem{paley2014bioluminescence}
{\sc M.~A. Paley and J.~A. Prescher}, {\em Bioluminescence: a versatile
  technique for imaging cellular and molecular features}, MedChemComm, 5
  (2014), pp.~255--267.

\bibitem{panda2019arrival}
{\sc S.~Panda}, {\em The arrival of circadian medicine}, Nature Reviews
  Endocrinology, 15 (2019), pp.~67--69.

\bibitem{ruben2019dosing}
{\sc M.~D. Ruben, D.~F. Smith, G.~A. FitzGerald, and J.~B. Hogenesch}, {\em
  Dosing time matters}, Science, 365 (2019), pp.~547--549.

\bibitem{sarkka2013bayesian}
{\sc S.~S{\"a}rkk{\"a}}, {\em Bayesian filtering and smoothing}, no.~3,
  Cambridge university press, 2013.

\bibitem{scheer2010impact}
{\sc F.~A. Scheer, K.~Hu, H.~Evoniuk, E.~E. Kelly, A.~Malhotra, M.~F. Hilton,
  and S.~A. Shea}, {\em Impact of the human circadian system, exercise, and
  their interaction on cardiovascular function}, Proceedings of the National
  Academy of Sciences, 107 (2010), pp.~20541--20546.

\bibitem{schnoerr2017approximation}
{\sc D.~Schnoerr, G.~Sanguinetti, and R.~Grima}, {\em Approximation and
  inference methods for stochastic biochemical kinetics—a tutorial review},
  Journal of Physics A: Mathematical and Theoretical, 50 (2017), p.~093001.

\bibitem{sethian1985curvature}
{\sc J.~A. Sethian}, {\em Curvature and the evolution of fronts},
  Communications in Mathematical Physics, 101 (1985), pp.~487--499.

\bibitem{stone2019application}
{\sc J.~E. Stone, X.~L. Aubert, H.~Maass, A.~J. Phillips, M.~Magee, M.~E.
  Howard, S.~W. Lockley, S.~M. Rajaratnam, and T.~L. Sletten}, {\em Application
  of a limit-cycle oscillator model for prediction of circadian phase in
  rotating night shift workers}, Scientific reports, 9 (2019), pp.~1--12.

\bibitem{stone2020computational}
{\sc J.~E. Stone, S.~Postnova, T.~L. Sletten, S.~M. Rajaratnam, and A.~J.
  Phillips}, {\em Computational approaches for individual circadian phase
  prediction in field settings}, Current Opinion in Systems Biology, 22 (2020),
  pp.~39--51.

\bibitem{sulli2018training}
{\sc G.~Sulli, E.~N. Manoogian, P.~R. Taub, and S.~Panda}, {\em Training the
  circadian clock, clocking the drugs, and drugging the clock to prevent,
  manage, and treat chronic diseases}, Trends in pharmacological sciences, 39
  (2018), pp.~812--827.

\bibitem{sun2008extended}
{\sc X.~Sun, L.~Jin, and M.~Xiong}, {\em Extended kalman filter for estimation
  of parameters in nonlinear state-space models of biochemical networks}, PloS
  one, 3 (2008), p.~e3758.

\bibitem{takahashi2017transcriptional}
{\sc J.~S. Takahashi}, {\em Transcriptional architecture of the mammalian
  circadian clock}, Nature Reviews Genetics, 18 (2017), pp.~164--179.

\bibitem{van1992stochastic}
{\sc N.~G. Van~Kampen}, {\em Stochastic processes in physics and chemistry},
  vol.~1, Elsevier, 1992.

\bibitem{walch2019sleep}
{\sc O.~Walch, Y.~Huang, D.~Forger, and C.~Goldstein}, {\em Sleep stage
  prediction with raw acceleration and photoplethysmography heart rate data
  derived from a consumer wearable device}, Sleep, 42 (2019), p.~zsz180.

\bibitem{walch2016global}
{\sc O.~J. Walch, A.~Cochran, and D.~B. Forger}, {\em A global quantification
  of “normal” sleep schedules using smartphone data}, Science advances, 2
  (2016), p.~e1501705.

\bibitem{wang2021level}
{\sc N.~Wang and D.~Forger}, {\em The level set kalman filter for state
  estimation of continuous-discrete systems}, IEEE Transactions on Signal
  Processing,  (2021).

\bibitem{wichniak2017treatment}
{\sc A.~Wichniak, K.~S. Jankowski, M.~Skalski, K.~Skwar{\l}o-So{\'n}ta, J.~B.
  Zawilska, M.~{\.Z}arowski, E.~Poradowska, and W.~Jernajczyk}, {\em Treatment
  guidelines for circadian rhythm sleep-wake disorders of the polish sleep
  research society and the section of biological psychiatry of the polish
  psychiatric association. part i. physiology, assessment and therapeutic
  methods}, Psychiatr Pol, 51 (2017), pp.~793--814.

\bibitem{woelders2017daily}
{\sc T.~Woelders, D.~G. Beersma, M.~C. Gordijn, R.~A. Hut, and E.~J. Wams},
  {\em Daily light exposure patterns reveal phase and period of the human
  circadian clock}, Journal of biological rhythms, 32 (2017), pp.~274--286.

\bibitem{zhang2014circadian}
{\sc R.~Zhang, N.~F. Lahens, H.~I. Ballance, M.~E. Hughes, and J.~B.
  Hogenesch}, {\em A circadian gene expression atlas in mammals: implications
  for biology and medicine}, Proceedings of the National Academy of Sciences,
  111 (2014), pp.~16219--16224.

\bibitem{zhu2012circadian}
{\sc L.~Zhu and P.~C. Zee}, {\em Circadian rhythm sleep disorders}, Neurologic
  clinics, 30 (2012), pp.~1167--1191.

\bibitem{zuurbier2015fragmentation}
{\sc L.~A. Zuurbier, A.~I. Luik, A.~Hofman, O.~H. Franco, E.~J. Van~Someren,
  and H.~Tiemeier}, {\em Fragmentation and stability of circadian activity
  rhythms predict mortality: the rotterdam study}, American journal of
  epidemiology, 181 (2015), pp.~54--63.

\end{thebibliography}
\end{document}